\documentclass[aip, cha, sd, amsmath,amssymb, preprint]{revtex4-1}

\newcommand  \eps{ \varepsilon}
\usepackage{graphicx}
\usepackage{dcolumn}
\usepackage{bm}

\begin{document}

\preprint{AIP/123-QED}

\title[]{Tangency bifurcation of invariant manifolds in a slow-fast system}

\author{Ian Lizarraga}
 \email{iml32@cornell.edu.}
\affiliation{Center for Applied Mathematics, Cornell University, Ithaca, NY 14853}

\date{\today}

\begin{abstract}
We study a three-dimensional dynamical system in two slow variables and one fast variable. We analyze the tangency of the unstable manifold of an equilibrium point with ``the'' repelling slow manifold, in the presence of a stable periodic orbit emerging from a Hopf bifurcation. This tangency heralds complicated and chaotic mixed-mode oscillations. We classify these solutions by studying returns to a two-dimensional cross section. We use the intersections of the slow manifolds as a basis for partitioning the section according to the number and type of turns made by trajectory segments. Transverse homoclinic orbits are among the invariant sets serving as a substrate of the dynamics on this cross-section. We then turn to a one-dimensional approximation of the global returns in the system, identifying saddle-node and period-doubling bifurcations. These are interpreted in the full system as bifurcations of mixed-mode oscillations. Finally, we contrast the dynamics of our one-dimensional approximation to classical results of the quadratic family of maps. We describe the transient trajectory of a critical point of the map over a range of parameter values.
\end{abstract}

\pacs{ 05.45.-a, 05.45.Ac}
\keywords{Singular Hopf bifurcation, tangency bifurcation of invariant manifolds, mixed-mode oscillations}
\maketitle

\begin{quotation}
We study a three-dimensional multiple timescale system in five parameters. A startling variety of behaviors can be identified as its five parameters are varied. Organizing this variety are the interactions between classical invariant manifolds (including fixed points, periodic orbits, and their (un)stable manifolds) and locally invariant slow manifolds. Here we focus on the interaction between the two-dimensional unstable manifold of a saddle-focus equilibrium point and a two-dimensional repelling slow manifold, in the presence of a stable periodic orbit of small amplitude.

 The images of global return maps, defined on carefully chosen two-dimensional cross-sections, are organized by the interactions of the attracting and repelling slow manifolds with these cross-sections. They are also influenced by the basin of attraction of the periodic orbit. We construct a symbolic map which partitions one such section according to the number and type of turning behaviors of the corresponding trajectories. We locate transverse homoclinic orbits to saddle points. On another cross-section, global returns are well-approximated by one-dimensional, nearly unimodal maps.  We show that saddle-node bifurcations of periodic orbits and period-doubling cascades occur.  Finally, we describe the dynamics of the critical point of the return map at carefully chosen parameters. 

Taking a broader view, our numerical results continue to point to the fruitful connections that exist between multiple-timescale flows and low-dimensional maps.
\end{quotation}

\section{\label{sec:intro} Introduction}

We study slow-fast dynamical systems of the form

\begin{eqnarray*}
\eps\dot{x} &=& f(x,y,\eps)\\
\dot{y} &=& g(x,y,\eps),
\end{eqnarray*}

where $x \in R^m$ is the {\it fast} variable, $y \in R^n$ is the {\it slow} variable, $\eps$ is the {\it singular perturbation parameter} that characterizes the ratio of the timescales,  and $f,g$ are sufficiently smooth. The {\it critical manifold} $C = \{f  = 0\}$ is the manifold of equilibria of the fast subsystem defined by $\dot{x} = f(x,y,0)$. When $\eps > 0$ is sufficiently small, theorems of Fenichel\cite{fenichel1972} guarantee the existence of locally invariant {\it slow manifolds} that perturb from subsets of $C$ where the equilibria are hyperbolic. We may also project the vector field $\dot{y} = g(x,y,0)$  onto the tangent bundle $TC$. Away from folds of $C$, we may desingularize this projected vector field to define the {\it slow flow}. The desingularized slow flow is oriented to agree with the full vector field near stable equilibria of $C$. For sufficiently small values of $\eps$, trajectories of the full system can be decomposed into segments lying on the slow manifolds near $C$ together with fast jumps across branches of $C$. Trajectory segments lying near the slow manifolds converge to solutions of the slow flow as $\eps$ tends to 0. 

\begin{figure*}
(a) \includegraphics[width=0.45\textwidth]{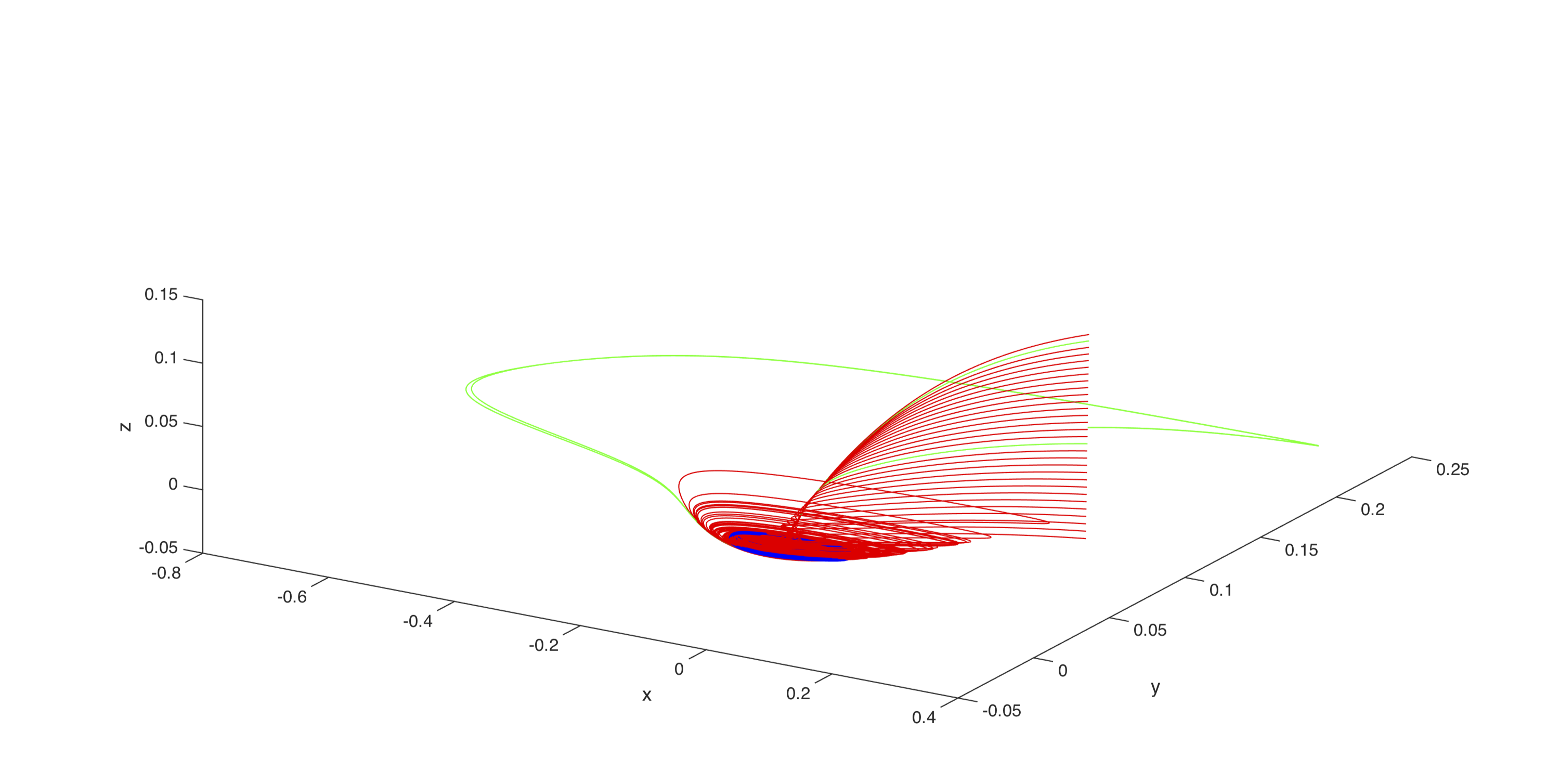}
(b) \includegraphics[width=0.45\textwidth]{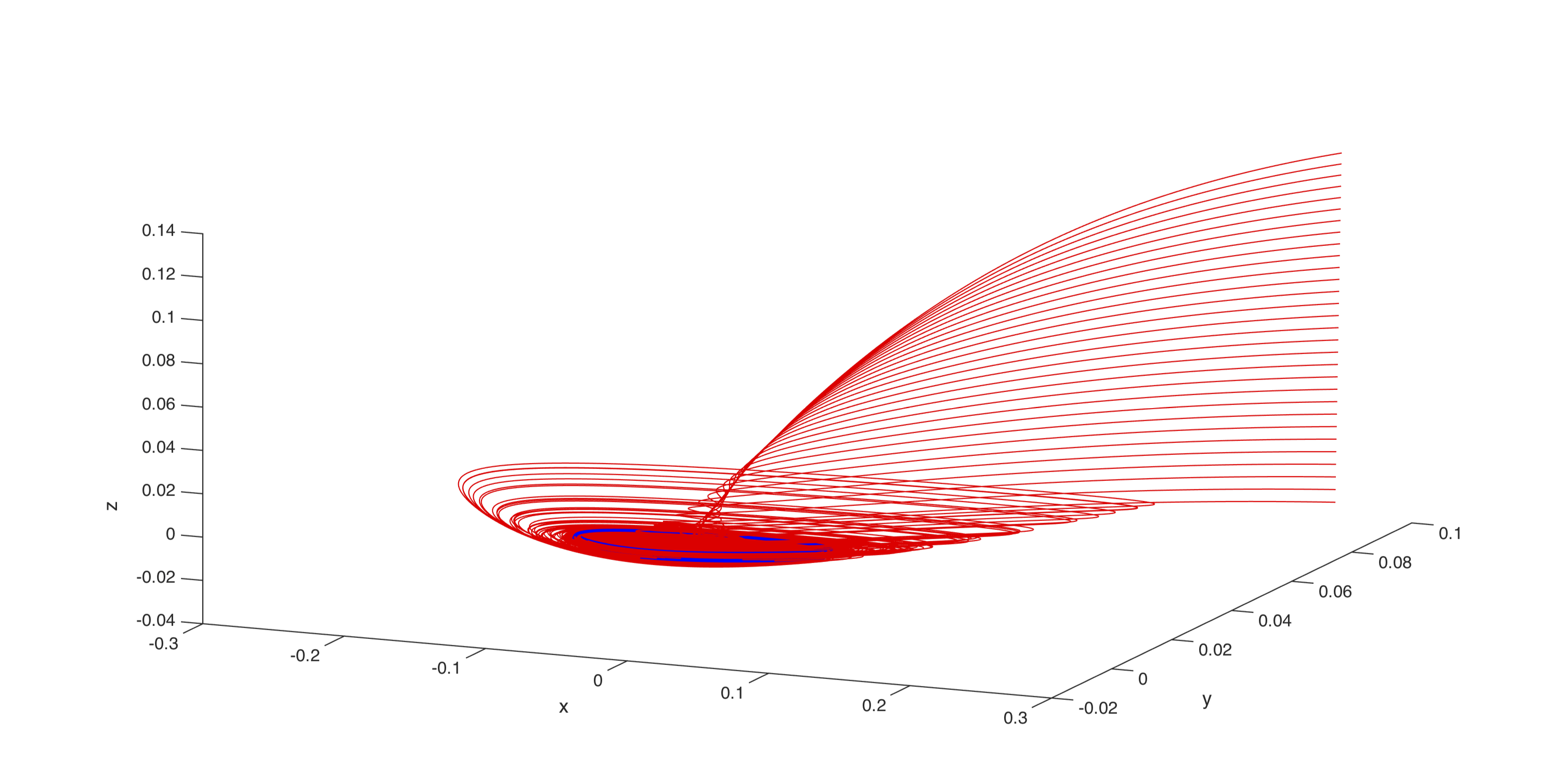}
\caption{\label{fig:tangbif} Phase space (a) just before ($\nu = 0.00647$) and (b) just after ($\nu = 0.00648$) tangency bifurcation of $W^u$ with $S^r_{\eps}$. Thirty trajectories are initialized in a band on $S^{a+}_{\eps}\cap \Sigma$, where $\Sigma = \{x = 0.27\}$. Blue curve: small-amplitude stable periodic orbit $\Gamma$. Red curves: forward trajectories tending asymptotically to $\Gamma$ without jumping to $S^{a-}_{\eps}$. Green curves: forward trajectories making a large-amplitude excursion before returning to $\Sigma$. Remaining parameters are $\eps = 0.01$, $a = -0.3$, $b = -1$, $c = 1$.}
\end{figure*}

We now focus on the case of two slow variables and one fast variable ($m=1$, $n=2$). The critical manifold $C$ is two-dimensional and folds of $C$ form curves. Points on fold curves are called {\it folded singularities}. when the slow flow is two-dimensional we use the terms ``folded node'', ``folded focus'', and ``folded saddle'' to denote folded singularities of node-, focus-, and saddle-type, respectively.  In analogy to classical bifurcation theory, folded saddle-nodes are folded singularities having a zero eigenvalue. When they exist, folded saddle-nodes are differentiated by whether they persist as equilibria in the full system of equations. We are interested here in folded-saddle nodes of type II (FSNII), which are true equilibria of the full system. It can be shown that {\it singular Hopf bifurcations} occur generically at distances $O(\eps)$ from the FSNII bifurcation in parameter space.\cite{guckenheimer2008siam} At this bifurcation, a pair of eigenvalues of the linearization of the flow crosses the imaginary axis, and a small-amplitude periodic orbit is born at the bifurcation point.

 Normal forms are used to study the local flow of full systems in neighborhoods of these folded singularities. Previous work by Guckenheimer\cite{guckenheimer2008chaos} analyzes the local flow maps and return maps of three-dimensional systems containing folded nodes and folded saddle-nodes. There, it is shown that the appearance of these folded singularities can give rise to complex and chaotic behavior. Characterizing the emergence of small-amplitude oscillations near a folded singularity has also been the subject of intense study. In the case of a folded node, Beno\^it\cite{benoit1990} and Wechselberger\cite{wechselberger2005} observed that the maximum number of small oscillations made by a trajectory passing through the folded node region is related to the ratio of eigenvalues of the folded node.

The present paper focuses on a dynamical system, defined in Sec. \ref{sec:shnf}, which contains folded singularities lying along a cubic critical manifold.  The critical manifold serves as a global return mechanism. Parametric subfamilies of this dynamical system have served as important prototypical models of electrochemical oscillations, including the Koper model\cite{koper1992}. This system serves as a concrete, minimal example of a three-dimensional system having an $S$-shaped critical manifold as a global return mechanism. Trajectories leaving a neighborhood of the folded singularities do so by jumping between branches of the critical manifold, before ultimately being reinjected into the regions containing the folded singularities. This interplay between local and global mechanisms gives rise to {\it mixed-mode oscillations} (MMOs), which are periodic solutions of the dynamical system containing large and small amplitudes and a distinct separation between the two. These solutions may be characterized by their signatures, which are symbolic sequences of the form $L_1^{s_1}L_2^{s_2} \cdots$. This notation is used to indicate that a particular solution undergoes $L_1$ large oscillations, followed by $s_1$ small oscillations, followed by $L_2$ large oscillations, and so on. The distinction between `large' and `small' oscillations is dependent on the model. Nontrivial aperiodic solutions are referred to as {\it chaotic MMOs}, and may be characterized as limits of families of MMOs as the lengths of the signatures grow very large.

The classification of routes to MMOs with complicated signatures as well as chaotic MMOs continues to garner interest. Global bifurcations have been identified as natural starting points in this direction. Even so, the connection between these bifurcations and  interactions of slow manifolds---which organize the global dynamics for small values of $\eps$---remains poorly understood. Period-doubling cascades, torus bifurcations,\cite{guckenheimer2008siam} and most recently, Shilnikov homoclinic bifurcations,\cite{guckenheimer2015} have been shown to produce MMOs with complex signatures. In the last case, one-dimensional approximations of return maps were used to analyze a Shilnikov bifurcation in a system which exhibits singular Hopf bifurcation.

In this paper, we use a similar technique to analyze a tangency of invariant manifolds.  Our starting point is a study by Guckenheimer and Meerkamp\cite{guckenheimer2012siam}, which comprehensively classifies local and global unfoldings of singular Hopf bifurcation. We describe the changes in the phase space as the unstable manifold of the saddle-focus equilibrium point crosses the repelling slow manifold of the system. Our approach takes for granted the complicated crossings of these two-dimensional manifolds, instead focusing directly on the influence of these crossings on the global dynamics. The main tool in our analysis is the approximation of the two-dimensional return map by a map on an interval, which parametrizes trajectories beginning on the attracting slow manifold. We show that in the presence of a small-amplitude stable periodic orbit, the one-dimensional return map has a rich topology. The domain of the map is disconnected, with components separated by finite-length gaps. Intervals where the return map is undefined correspond to bands of initial conditions in the full system whose forward trajectories asymptotically approach the small-amplitude stable periodic orbit without making a large-amplitude passage. The first and second derivatives of the map grow very large outside of large subintervals where the map is unimodal. 

We also interpret classical bifurcations of the one-dimensional map as routes to chaotic behavior in the full system. We show that a period-doubling cascade occurs in this map, which gives rise to chaotic MMOs. This cascade is reminiscent of the classical cascade in the family of quadratic maps, even though on small subsets, our return map is far from unimodal. Saddle-nodes of mixed-mode cycles, defined as fixed points of the return map with unit derivative, are also shown to occur. Finally, we identify a parameter set for which the full dynamics is close to the dynamics of a unimodal map with a critical point having dense forward orbit.

\section{\label{sec:shnf} Three-Dimensional System of Equations}

We study the following three-dimensional flow:
\begin{eqnarray}
\eps \dot{x} &=& y - x^2 - x^3 \nonumber\\
\dot{y} &=& z - x \label{eq:shnf}\\
\dot{z} &=& -\nu - ax -by - cz,\nonumber
\end{eqnarray}

where $x$ is the fast variable, $y,z$ are the slow variables, and $\eps,\nu,a,b,c$ are the system parameters. This system exhibits a singular Hopf bifurcation.\cite{braaksma1998,guckenheimer2008siam,guckenheimer2012dcds} The critical manifold is the S-shaped cubic surface $C = \{y = x^2 + x^3\}$ having two fold lines $L_0 := S \cap \{x = 0\}$ and $L_{-2/3} := S \cap \{ x = -2/3\}$. When $\eps > 0 $ is sufficiently small, nonsingular portions of $C$ perturb to families of slow manifolds: near the branches $S\cap \{x > 0\}$ (resp. $S \cap \{x < -2/3\}$), we obtain the {\it attracting slow manifolds} $S^{a+}_{\eps}$ (resp. $S^{a-}_{\eps}$) and near the branch $S\cap\{-2/3<x<0\}$ we obtain the {\it repelling slow manifolds} $S^r_{\eps}$. Nearby trajectories are exponentially attracted toward $S^{a\pm}_{\eps}$ and exponentially repelled from $S^r_{\eps}$. One derivation of these estimates uses the Fenichel normal form.\cite{jones1994} Within each family, these sheets are $O(-\exp(c/ \eps))$ close\cite{jones1994,jones1995}, so we refer to any member of a particular family as `the' slow manifold. This convention should not cause confusion. 

We focus on parameters where forward trajectories beginning on $S^{a+}_{\eps}$ interact with a `twist region' near $L_0$, a saddle-focus equilibrium point $p_{eq}$, or both. A folded singularity $n =(0,0,0) \in L_0$ is the governing center of this twist region. The saddle-focus $p_{eq}$ has a two-dimensional unstable manifold $W^u$ and a one-dimensional stable manifold $W^s$. This notation disguises the dependence of these manifolds on the parameters of the system.

\section{Tangency bifurcation of invariant manifolds}

Guckenheimer and Meerkamp\cite{guckenheimer2012siam} drew bifurcation diagrams of the system \eqref{eq:shnf} in a two-dimensional slice of the parameter space defined by $\eps = 0.01$, $b = -1$, and $c = 1$. Codimension-one tangencies of $S^r_{\eps}$ and $W^u$ are represented in Figure 5.1 of their paper by smooth curves (labeled T) in $(\nu,a)$ space. For fixed $a$ and increasing $\nu$, this tangency occurs after $p_{eq}$ undergoes a supercritical Hopf bifurcation. A parametric family of stable limit cycles emerges from this bifurcation. Henceforth we refer to `the' small-amplitude stable periodic orbit $\Gamma$ to refer to the corresponding member of this family at a particular parameter set. The two-dimensional stable manifolds of $\Gamma$ interact with the other invariant manifolds of the system. Guckenheimer and Meerkamp identify a branch of period-doubling bifurcations as $\nu$ continues to increase after the first slow-manifold tangency. We show that the basin of attraction of the periodic orbit has a significant influence on the global returns of the system.

Fixing $a = -0.03$, the tangency occurs within the range $\nu \in \left[ 0.00647, 0.00648\right]$. The location of the tangency may be approximated by studying the asymptotics of orbits beginning high up on $S^{a+}_{\eps}$. Fix a section $\Sigma = S^{a+}_{\eps}\cap \{x = 0.27\}$.  Before the tangency occurs, trajectories lying on and sufficiently near $W^u$ must either escape to infinity or asymptotically approach $\Gamma$; these trajectories cannot jump to the attracting branches of the slow manifold, as they must first intersect $S^r_{\eps}$ before doing so. Trajectories beginning in $\Sigma$ first flow very close to $p_{eq}$. As shown in Figure 1, these trajectories then leave the region close to $W^u$. We observe that before the tangency, $W^u$ forms a boundary of the basin of attraction of $\Gamma$. Therefore, all trajectories sufficiently high up on $S^{a+}_{\eps}$ must lie inside the basin of attraction (Figure \ref{fig:tangbif}a).

After the tangency has occurred, isolated trajectories lying in $W^u$ will also lie in $S^r_{\eps}$. These trajectories will bound sectors of trajectories which can now make large-amplitude passages. Trajectories within these sectors jump `to the left' toward $S^{a-}_{\eps}$ or `to the right' toward $S^{a+}_{\eps}$. Trajectories initialized in $\Sigma$ that leave neighborhoods of $p_{eq}$ near these sectors contain {\it canard} segments, which are solution segments lying along $S^r_{\eps}$. Examples of such trajectories are highlighted in green in Figure 1. We can now establish a dichotomy between those trajectories in $\Sigma$ that immediately flow to $\Gamma$ and never leave a small neighborhood of the periodic orbit, versus those that make a global return. In Figure \ref{fig:tangbif}b, only two of the thirty sample trajectories are able to make a global return. Near the boundaries of these subsets, trajectories can come arbitrarily close to $\Gamma$ before escaping and making one large return. Note however that such trajectories might still lie inside the basin of attraction of $\Gamma$, depending on where they return on $\Sigma$. Such trajectories escape via large-amplitude excursions at most finitely many times before tending asymptotically to $\Gamma$.  We now focus on the parameter regime where the tangency has already occurred. In Figure 5.1 of the paper of Guckenheimer and Meerkamp, this corresponds to the region to the right of the $T$ (manifold tangency) curve.

\begin{figure}
\includegraphics[width=0.5\textwidth]{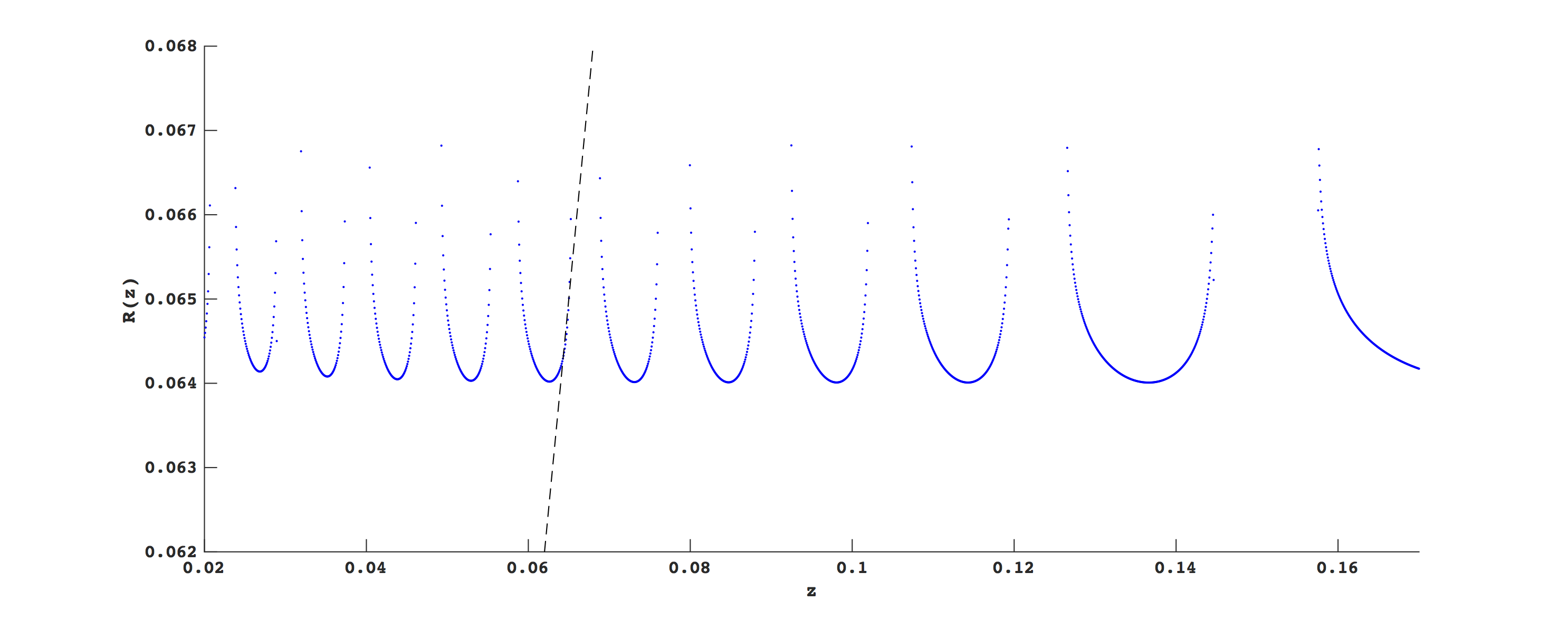}
\caption{\label{fig:retmap}   The return map $R: \Sigma_+ \to \Sigma_+$ of the system \eqref{eq:shnf} with $\Sigma_+ = \{x = 0.3\}$. Points in the two-dimensional section are parametrized by their $z$-coordinates. The dashed black line is the line of fixed points $\{(z,z)\}$. Parameter set: $\nu = 0.00802$, $a = -0.3$, $b = -1$, $c = 1$.}
\end{figure}

\section{\label{sec:maps} Singular and Regular Returns}

\begin{figure}
(a) \includegraphics[width=0.45\textwidth]{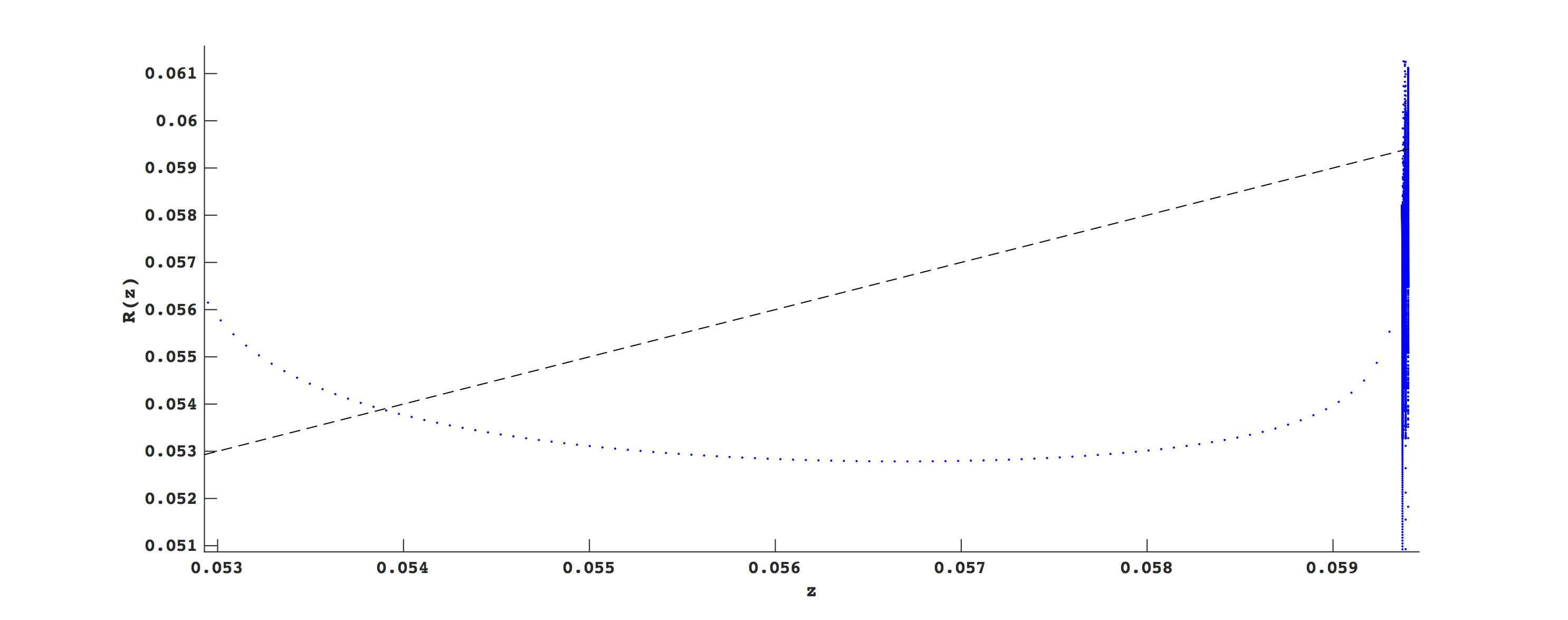}
(b) \includegraphics[width=0.45\textwidth]{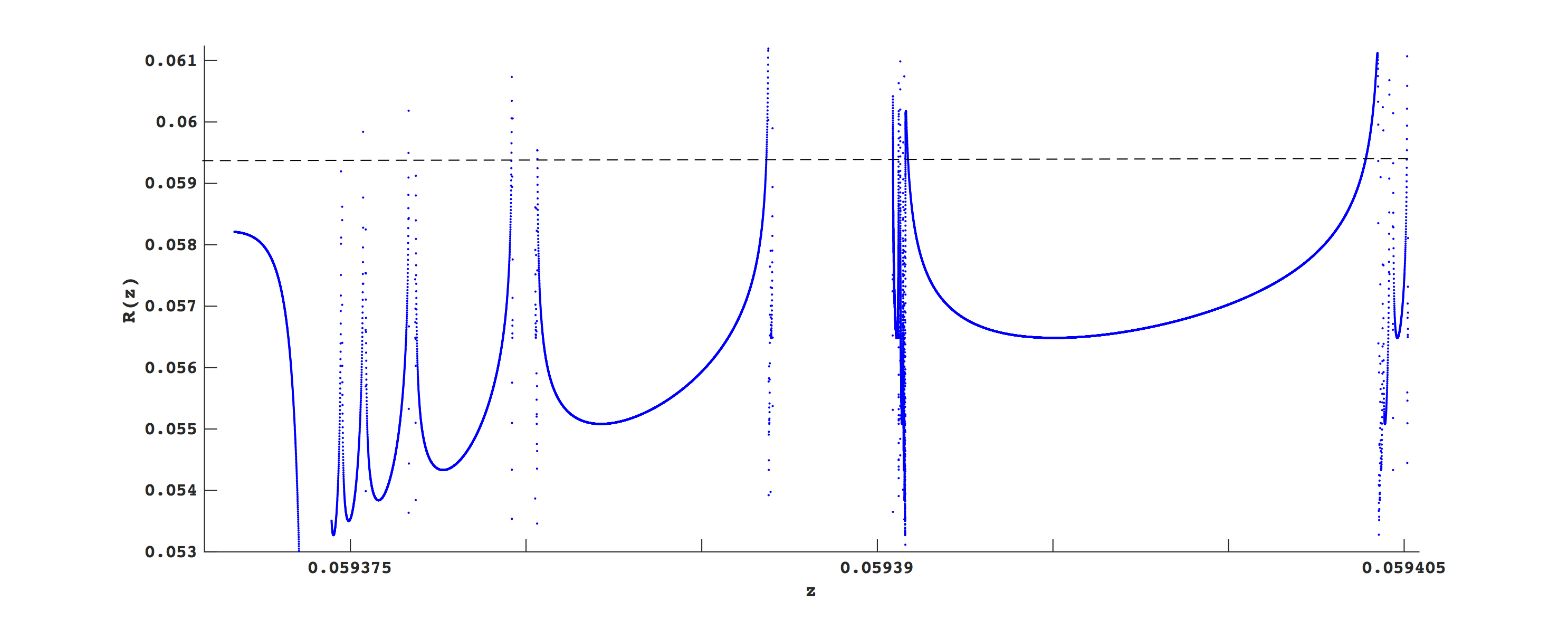}\\
(c) \includegraphics[width=0.45\textwidth]{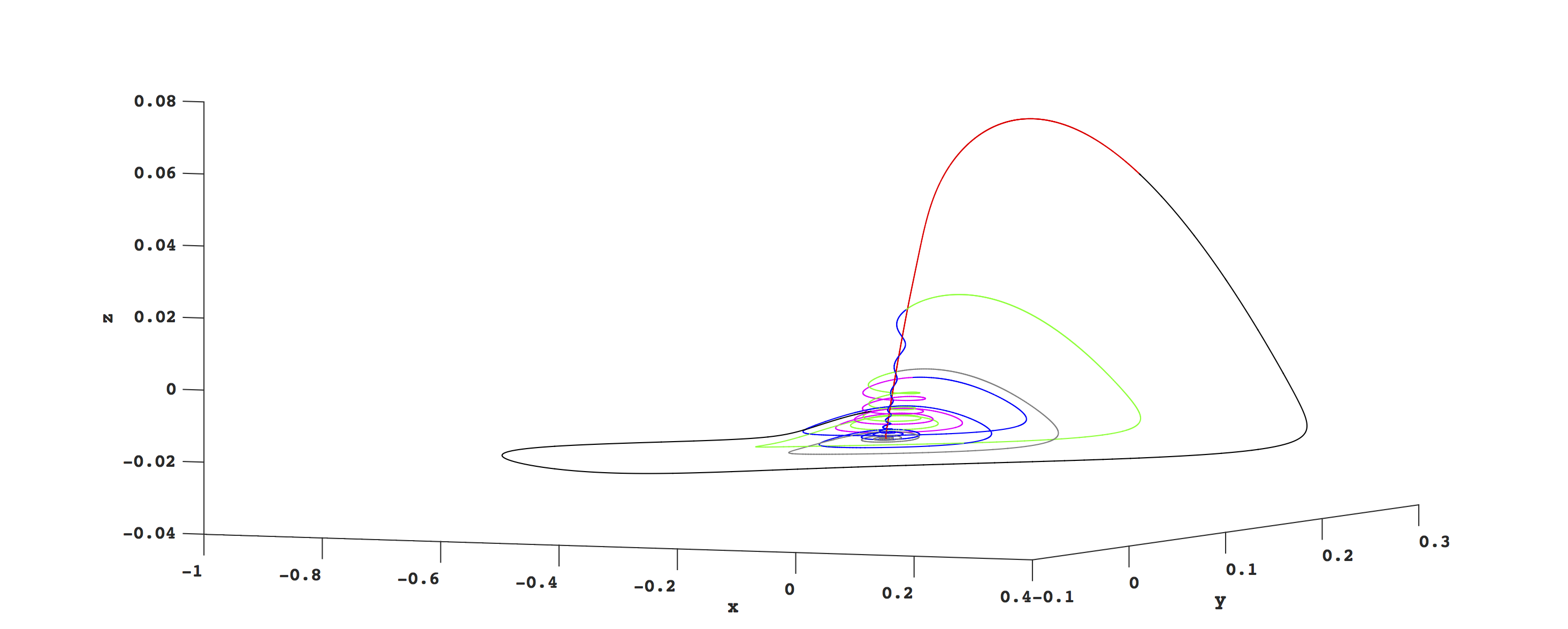}
(d) \includegraphics[width=0.45\textwidth]{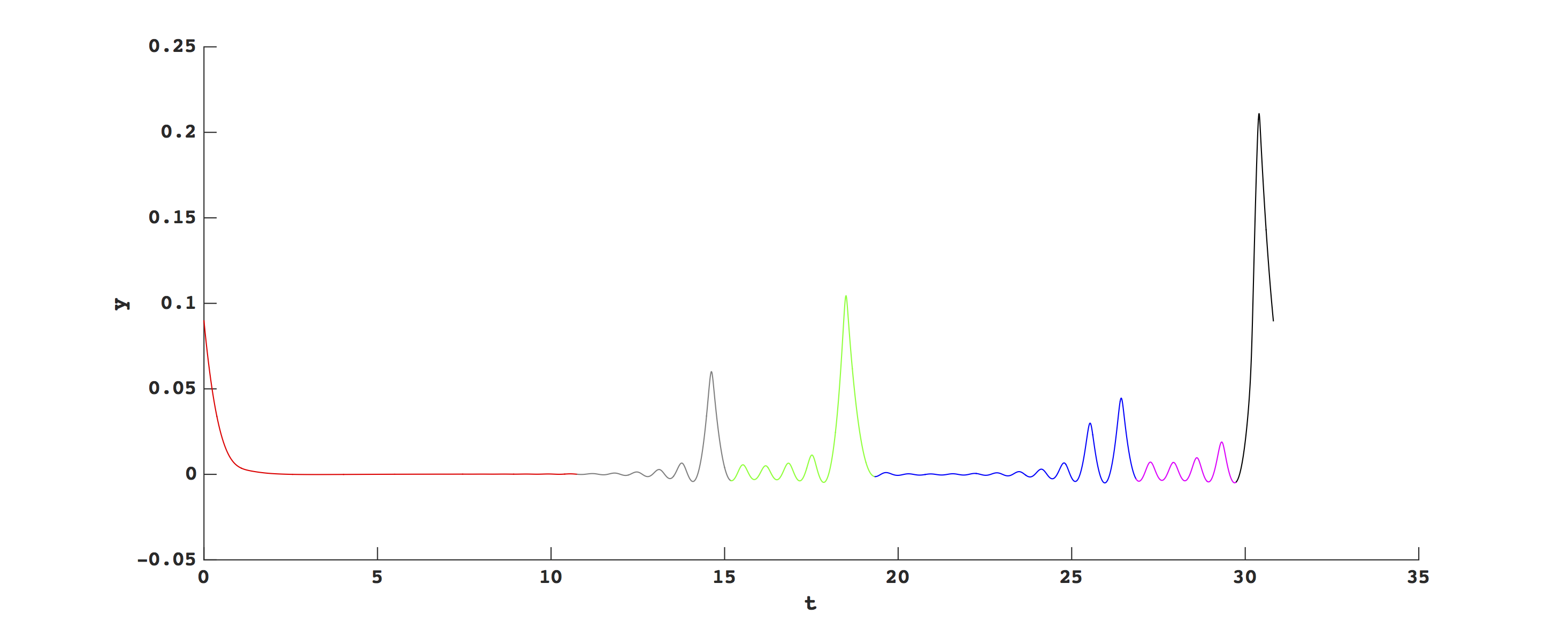}
\caption{\label{fig:retmap2}  (a) Subinterval of the return map $R: \Sigma_+ \to \Sigma)+$  of Eqs. \eqref{eq:shnf} and (b) a refinement of the subinterval. Dashed black line is the line of fixed points $R(z) = z$. (c) Periodic orbit corresponding to fixed point of $R$ at $z \approx 0.05939079$. (d) Time series of the periodic orbit. The orbit is decomposed into red, gray, green, blue, magenta, and black segments (defined as in Sec. \ref{sec:maps}). Parameter set: $\nu \approx 0.00870134$, $a = -0.3$, $b = -1$, $c = 1$.}
\end{figure}

Approximating points on $\Sigma$ by their $z$-coordinates, the return map $R: \Sigma \to \Sigma$ is well-approximated by a one-dimensional map on an interval, also denoted $R$. In the presence of the small-amplitude stable periodic orbit $\Gamma$, we now compare our one-dimensional approximation to return maps in the case of folded nodes\cite{wechselberger2005} and folded saddle-nodes\cite{guckenheimer2008chaos,krupa2010}. Where the return map is defined, trajectories beginning in different components of the domain of $R$ make different numbers of small turns before escaping the local region. These subsets are somewhat analogous to the rotation sectors arising from twists due to a folded node.\cite{wechselberger2005} However, in the present case there is a folded singularity as well as a saddle-focus as well as a small-amplitude periodic orbit. Each of these local objects plays a role in the twisting of trajectories that enter neighborhoods of the fold curve $L_0$.

When the small-amplitude stable periodic orbit exists, the domain of the return map is now disconnected, with components separated by finite-length gaps (Figure \ref{fig:retmap}). The gaps where $R$ is undefined correspond to those trajectories beginning on $S^{a+}_{\eps}$ that asymptotically approach $\Gamma$ without making a large-amplitude oscillation.  The second difference concerns the extreme nonlinearity near the boundaries of the disconnected intervals where $R$ is defined (Figure \ref{fig:retmap2}a). Portions of the image lie below the local minima in these local concave segments, resulting in tiny regions near the boundaries where the derivative changes rapidly. These points arise from canard segments of trajectories resulting in a jump from $S^r_{\eps}$ to $S^{a+}_{\eps}$ and hence to $\Sigma$.  Fixing the parameters and iteratively refining successively smaller intervals of initial conditions, this pattern of disconnected regions where the derivative changes rapidly seems to repeat up to machine accuracy. One consequence of this structure is the existence of large numbers of unstable periodic orbits, defined by fixed points of $R$ at which $|R'(z)| > 1$. This topological structure also appears to be robust to variations of the parameter $\nu$. 

This complicated structure arises from the interaction between the basin of attraction of $\Gamma$, the twist region near the folded singularity and $W^{u,s}$. As an illustration of this complexity, consider an unstable fixed point $z \approx  0.05939079$ of the return map as defined in Figure \ref{fig:retmap2}(b), interpreted as an unstable periodic orbit in the full system of equations (Figure \ref{fig:retmap2}(c)-(d)). The orbit is approximately decomposed according to its interactions with the (un)stable manifolds of $p_{eq}$ and the slow manifolds. One possible forward-time decomposition of this orbit proceeds as follows:

\begin{itemize}
\item A segment (red) that begins on $S^{a+}_{\eps}$ and flows very close to $p_{eq}$ by remaining near $W^s$,
\item a segment (gray) that leaves the region near $p_{eq}$ along $W^u$, then jumping right from $S^r$ to $S^{a+}_{\eps}$, 
\item a segment (green) that flows from $S^{a+}_{\eps}$ to $S^{r}_{\eps}$, making small-amplitude oscillations while remaining a bounded distance away from $p_{eq}$, then jumping right from $S^r_{\eps}$ to $S^{a+}_{\eps}$,
\item a segment (blue) that flows back down into the region near $p_{eq}$, making small oscillations around $W^s$, then jumping right from $S^r_{\eps}$ to $S^{a+}_{\eps}$,
\item a segment (magenta) with similar dynamics to the green segment, making small-amplitude oscillations while remaining a bounded distance away from $p_{eq}$, then jumping right from $S^r_{\eps}$ to $S^{a+}_{\eps}$, and
\item a segment (black) making a large-amplitude excursion by jumping left to $S^{a-}_{\eps}$, flowing to the fold $L_{-2/3}$, and then jumping to $S^{a+}_{\eps}$.
\end{itemize}

A linearized flow map can be constructed\cite{glendinning1984,silnikov1965} in small neighborhoods of the saddle-focus $p_{eq}$, which can be used to count the number of small-amplitude oscillations contributed by orbit segments approaching the equilibrium point. However, the small-amplitude periodic orbit and the twist region produce additional twists, as observed in the green and magenta segments of the example above. 

\begin{figure}
(a)\includegraphics[width=0.48\textwidth]{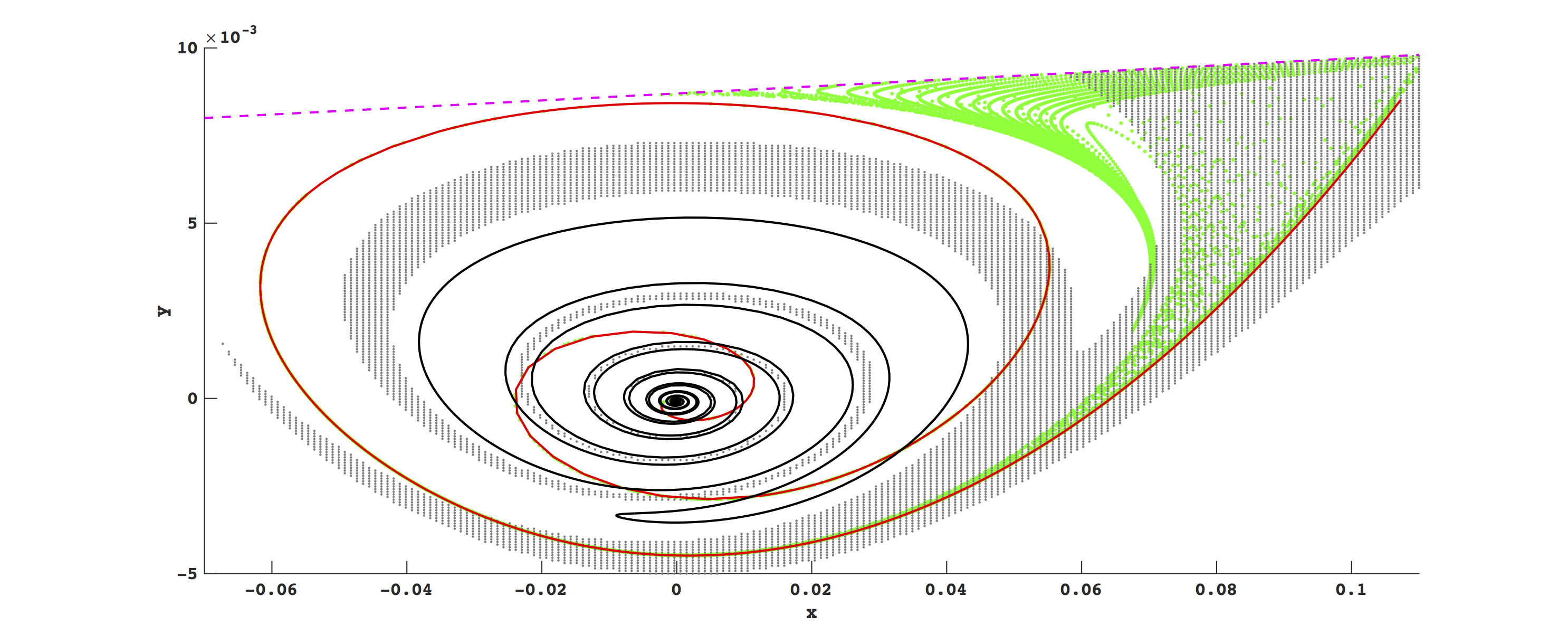}\\
(b)\includegraphics[width=0.48\textwidth]{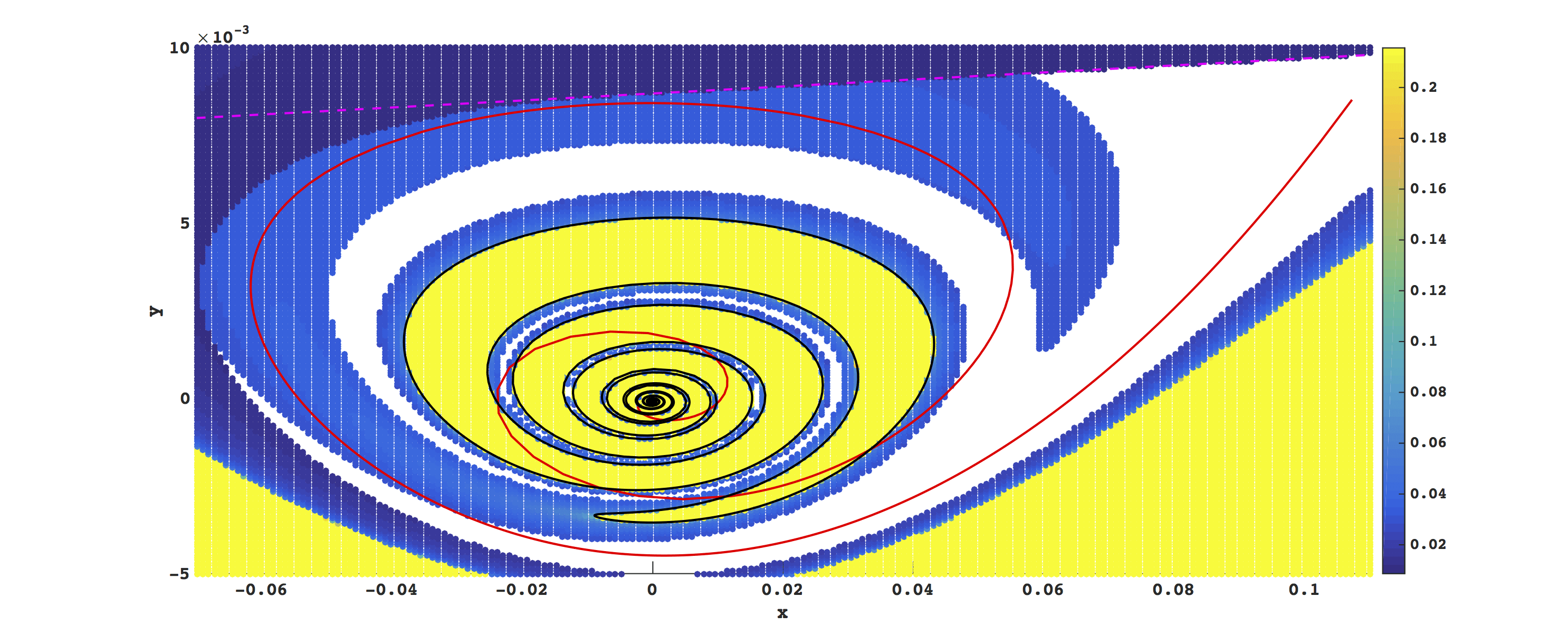}
\caption{\label{fig:z0} (a) Geometry in the section $\Sigma_0 = \{(x,y,z): x \in \left[-0.07,0.11\right], y \in \left[-0.005,0.01\right], z=0\}$.  Gray points sample the subset of $\Sigma_0$ whose corresponding forward trajectories tend asymptotically close to the stable periodic orbit without returning to $\Sigma_0$. Green points denote the first forward return of the remaining points in $\Sigma_0$ with the orientation $\dot{z} < 0$. (b) Color plot of maximal height ($y$-coordinate) obtained by trajectories that return to $\Sigma_0$ as defined in (a). Cross-sections of $S^{a+}_{\eps}$ (red) and $S^r_{\eps}$ (black) at $\Sigma_0$ are shown, and the tangency of the vector field with $\Sigma_0$ (i.e. the set $\{ax + by = -\nu\}$) is given by the magenta dashed line. Parameter set: $\nu \approx 0.00870134$, $a = 0.01$, $b = -1$, $c = 1$.}
\end{figure}

We will return to one-dimensional approximations of the return map in Sec. \ref{sec:ret}, but now we focus on two-dimensional maps, and show that we can illuminate key features of their small-amplitude oscillations. Let us fix a cross-section and define the geometric objects whose interactions organize the return dynamics. Define $\Sigma_0$ to be a compact subset of $\{ z = 0\}$ containing the first intersection (with orientation $\dot{z} >0$) of $W^s$ . Let $B_0$ denote the {\it immediate basin of attraction} of the stable periodic orbit $\Gamma$, which we define as the set of points in $\Sigma_0$ whose forward trajectories under the flow of Eq. \eqref{eq:shnf} asymptotically approach $\Gamma$ without returning to $\Sigma_0$, and let $\partial B_0$ denote its boundary. The periodic orbit itself does not intersect our choice of cross-section. 

Since we wish to study trajectory segments that return to the cross-section, $B_0$ functions as an escape subset. Rigorously, the forward return map $R: \Sigma_0 \to \Sigma_0$ is undefined on the subset $B_0$, and points landing in $B_0$ under forward iterates of $R$ `escape'.  Obviously the trajectories with initial conditions inside $\cup_{i=0}^{\infty}R^{-i}(B_0)$ are contained within the basin of attraction of $\Gamma$, and furthermore the $j$-th iterate of the return map $R^j$ is defined only on the subset $\Sigma_0 - \cup_{i=0}^j R^{-j}(B_0)$. Finally, we abuse notation slightly and denote by $S^{a+}_{\eps}$ (resp. $S^r_{\eps}$) the intersections of the corresponding slow manifolds with $\Sigma_0$. We also refer to the intersection of $S^{a+}_{\eps}$ (resp. $S^r_{\eps}$) with $\Sigma_0$ as the {\it attracting} (resp. {\it repelling}) {\it spiral} due to its distinctive shape (see Figure \ref{fig:z0}). The immediate basin of attraction $B_0$ is depicted by gray points in Figure \ref{fig:z0}(a). This result implies that the basin of attraction of the periodic orbit contains at least a thickened spiral which $S^{a+}_{\eps}$ intersects transversely in interval segments, accounting for the disconnected images of the one-dimensional return maps. 

The slow manifolds also intersect transversely. Segments of the attracting spiral can straddle both $B_0$ and the repelling spiral. In Fig. \ref{fig:z0}(b), we color initial conditions based on the maximum $y$-coordinate achieved by the corresponding trajectory before its return to $\Sigma_0$.  Due to the Exchange Lemma, only thin bands of trajectories are able to remain close enough to $S^r_{\eps}$ to jump at an intermediate height. 
We choose the maximum value of the $y$-coordinate to approximately parametrize the length of the canards. This parametrization heavily favors trajectories jumping left (from $S^r_{\eps}$ to $S^{a-}_{\eps}$) rather than right (from $S^r_{\eps}$ to $S^{a+}_{\eps}$), since trajectories jumping left can only return to $\Sigma_0$ by first following $S^{a-}_{\eps}$ to a maximal height, and then jumping from $L_{-2/3}$ to $S^{a+}_{\eps}$. This asymmetry is useful: in Figure \ref{fig:z0}, $S^r_{\eps}$ serves as a boundary between the (apparently discontinuous) blue and yellow regions, clearly demarcating those trajectories which turn right rather than left before returning to $\Sigma_0$. Summarizing, $\partial B_0$ and $S^r_{\eps}$ partition this section according to the behavior of orbits containing canards. 

\begin{figure}
(a) \includegraphics[width=0.45\textwidth]{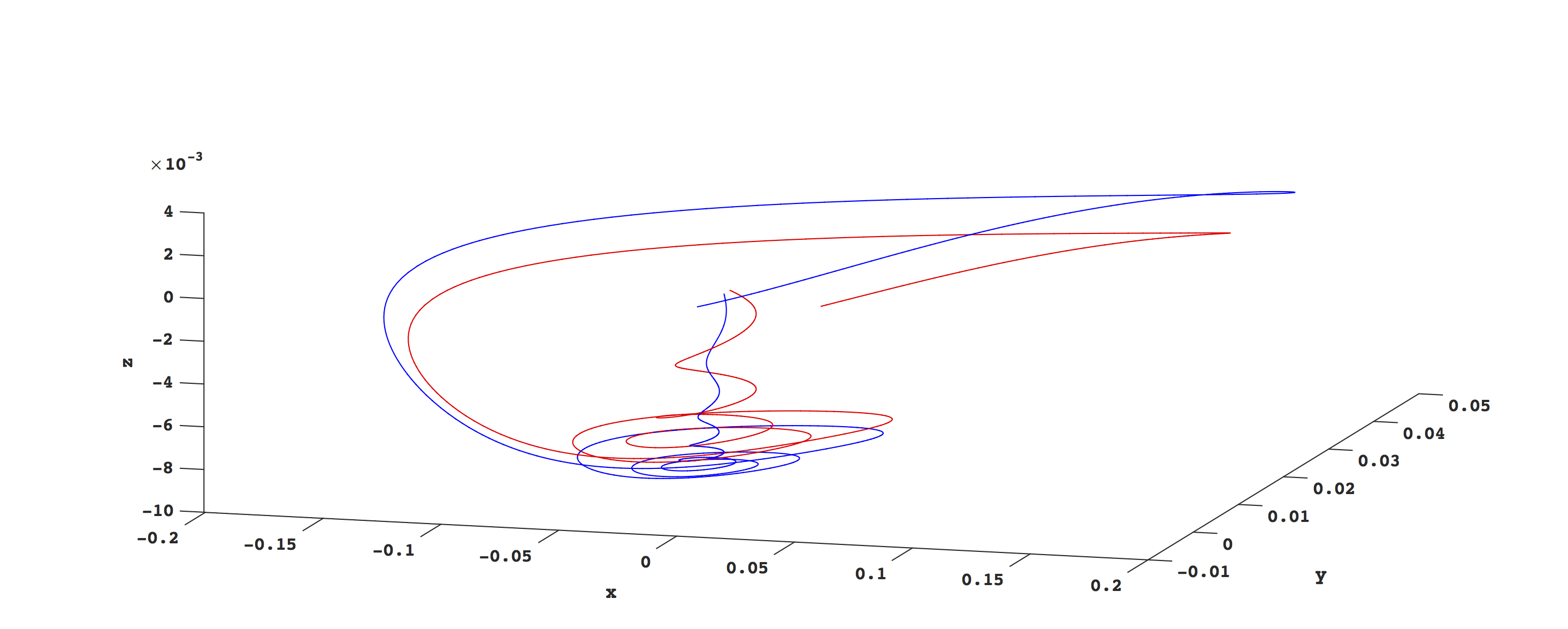}\\
(b) \includegraphics[width=0.45\textwidth]{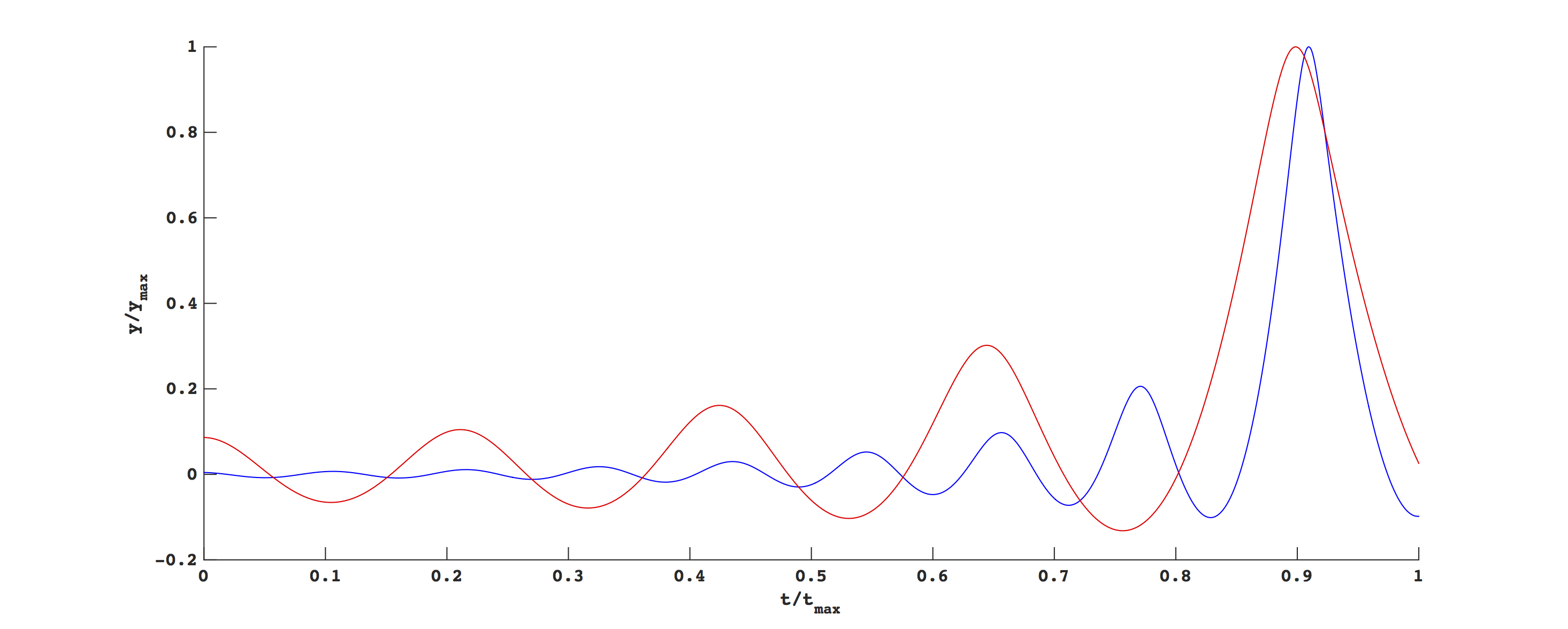}
\caption{\label{fig:sao}  (a) Two phase space trajectories beginning and ending on the section $\{z = 0\}$ with stopping condition $\dot{z} < 0$ and (b) the normalized time series of their normalized $y$-coordinates of each trajectory. Initial conditions: blue, $(x,y,z) = (0.000553, 0.000201, 0)$; red, $(x,y,z) = (0.000553, 0.003065, 0)$. Parameter set: $\nu \approx 0.00870134$, $a = 0.01$, $b = -1$, $c = 1$.}
\end{figure}

Trajectories beginning in $\Sigma_0$ either follow $W^s$ closely and spiral out along $W^u$ or remain a bounded distance away from both the equilibrium point and $W^s$, instead making small-amplitude oscillations consistent with a folded node. Differences between these two types of small-amplitude oscillations have been observed in earlier work. The transition from the first kind of small-amplitude oscillation to the second is a function of the distance from the initial condition to the intersection of $W^s$ with the cross-section. Two initial conditions are chosen on a vertical line embedded in the section $\{z=0\}$, having the property that the resulting trajectory jumps right from $S^r_{\eps}$ at an intermediate height before returning to the section with orientation $\dot{z} < 0$. These initial conditions are found by selecting points in Figure \ref{fig:z0}(b) in the blue regions lying on a ray that extends outward from the center of the repelling spiral. The corresponding return trajectories are plotted in Figure \ref{fig:sao}. The production of small-amplitude oscillations is dominated by the saddle-focus mechanism: in the example shown, the red orbit exhibits four oscillations before the (relatively) large-amplitude return, whereas the blue orbit exhibits seven oscillations. We can select trajectories with increasing numbers of small-amplitude oscillations by picking points closer to $W^s \cap \{z=0\}$.  A complication in this analysis is that jumps at intermediate heights, which are clearly shown to occur in these examples, blur the distinction between `large' and `small' oscillations in a mixed-mode cycle.

\begin{figure}
(a) \includegraphics[width=0.4\textwidth]{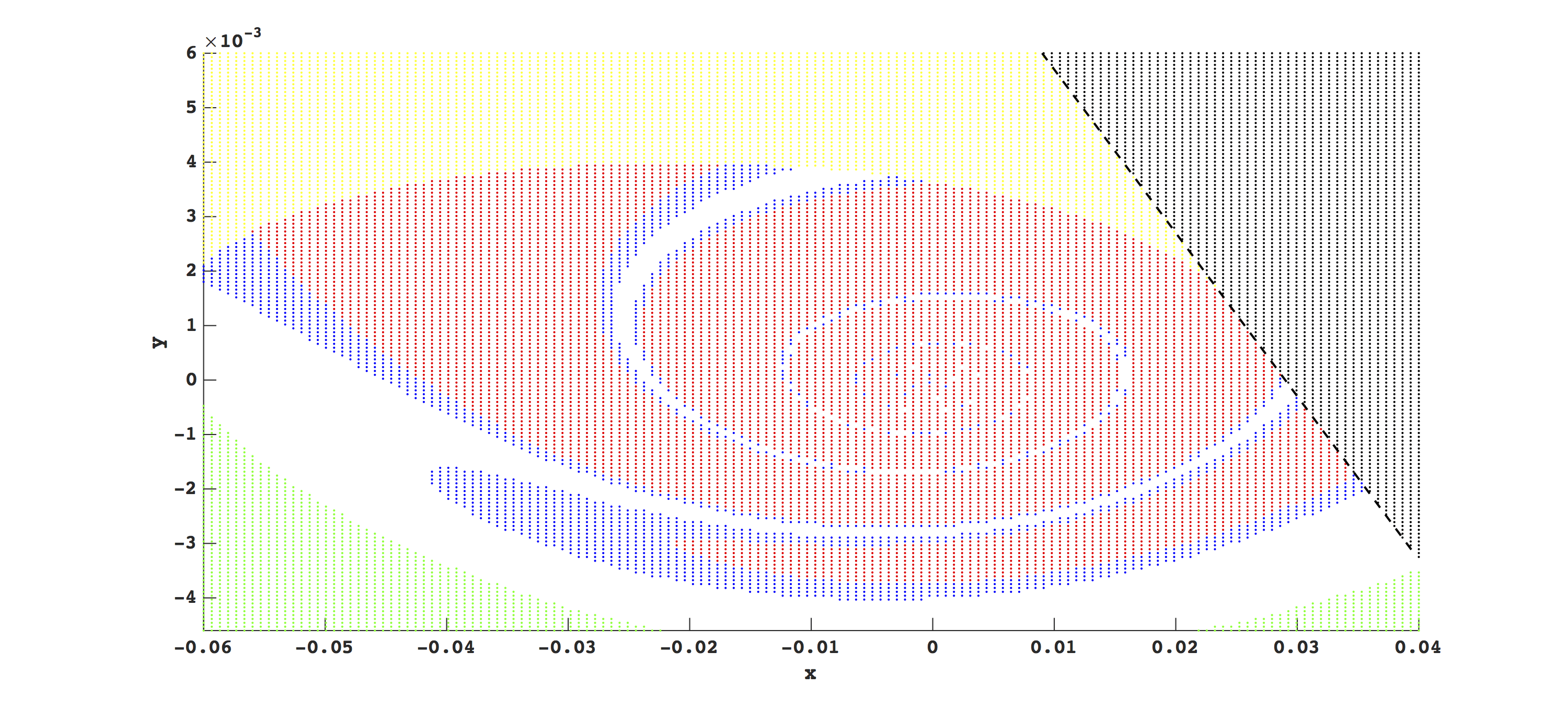}
(b) \includegraphics[width=0.4\textwidth]{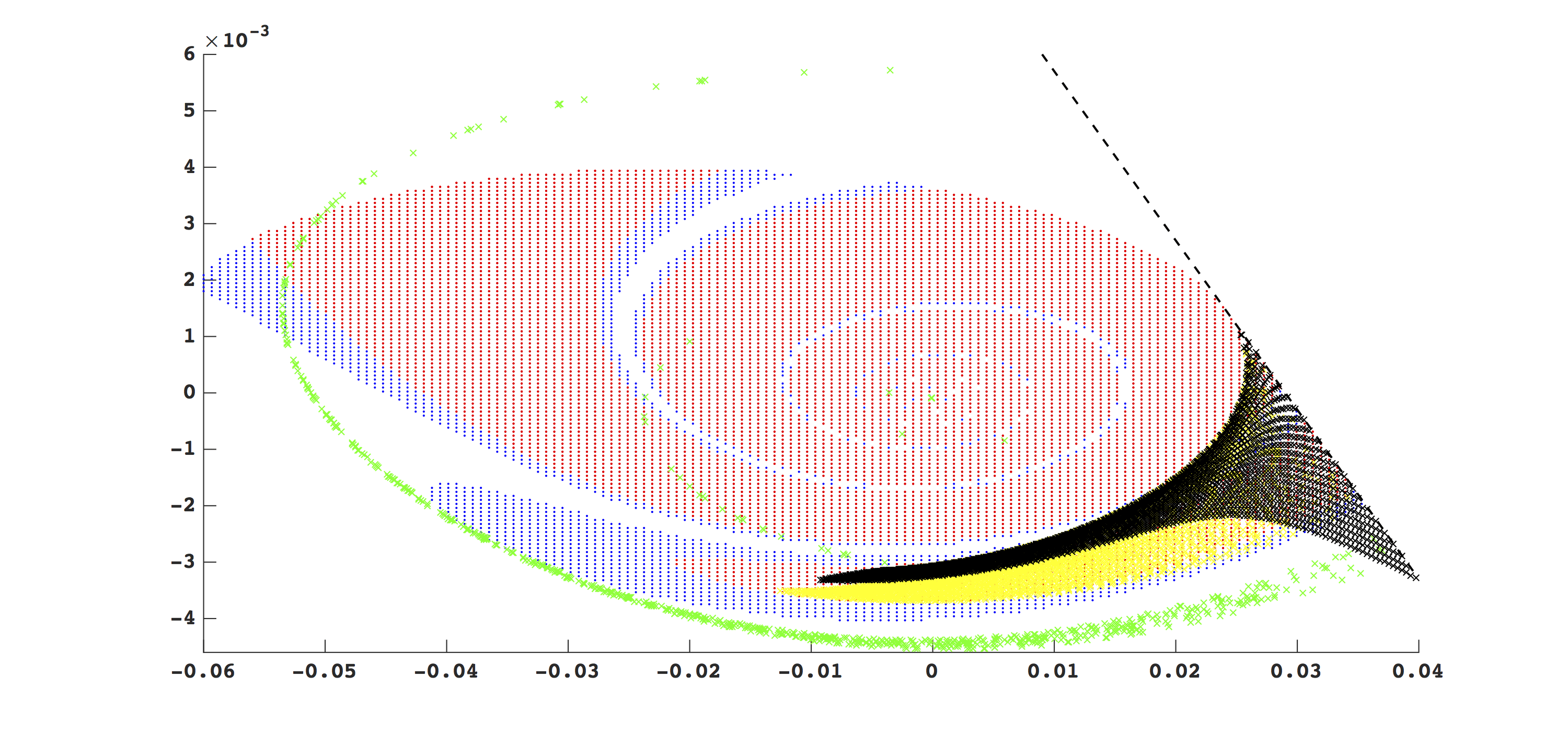}\\
(c) \includegraphics[width=0.4\textwidth]{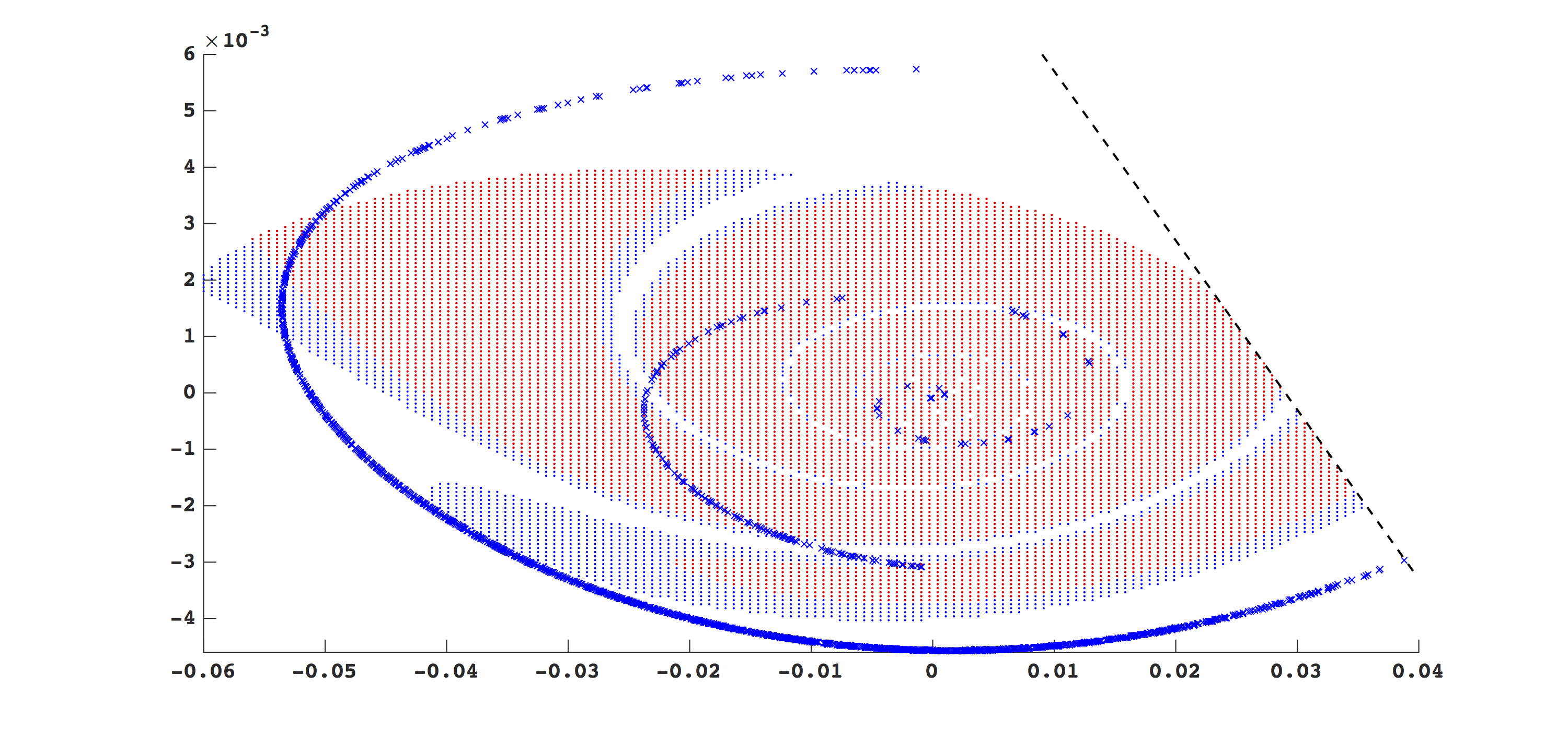}
(d) \includegraphics[width=0.4\textwidth]{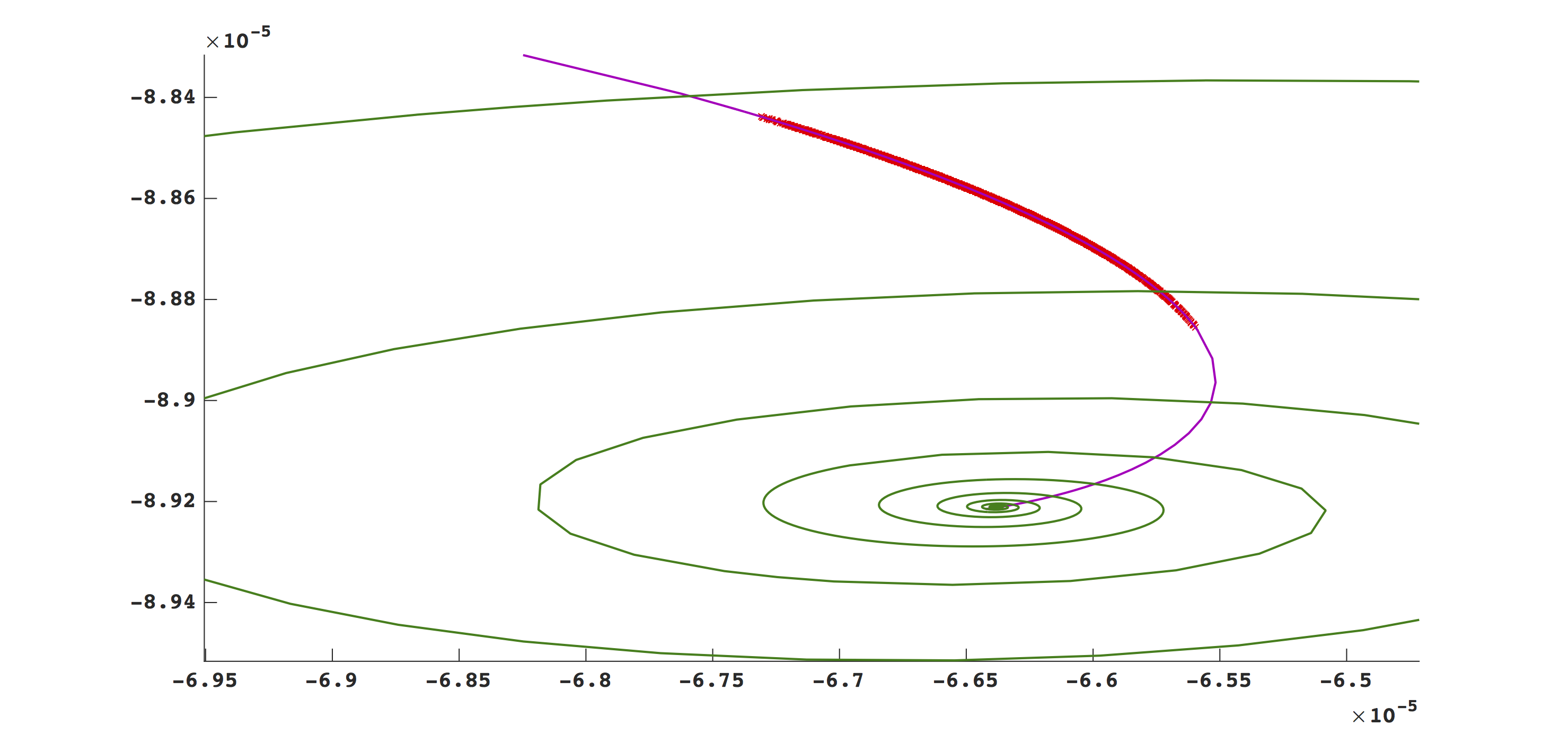}
\caption{\label{fig:2dmap} (a) Partition of a compact subset of the cross-section $\Sigma_0 = \{z = 0\}$. Black dashed line is the tangency of the vector field $\{\dot{z} = 0\}$, separating the subset $\{\dot{z} > 0\}$ (black points) from the other subsets. Yellow (resp. green): points above (resp. below) the line $\{y = 0\}$ with winding number less than three. Red (resp. blue): points whose forward trajectories reach a maximal height greater than (resp. less than) 0.18 and have winding number three or greater. (b) Overlay of red and blue subsets of domain (points) with images of yellow, green, and black subsets (crosses). (c) Overlay of red and blue subsets of domain (points) with the image of the blue subset (crosses). (d) Overlay of attracting spiral (magenta), repelling spiral (dark green), and image of red subset (crosses). Note the change in scale of the final figure. Generated from a $500 \times 500$ grid of initial conditions beginning on $\Sigma_0$. Parameter set: $\nu \approx 0.00870134$, $a = -0.3$, $b = -1$, $c = 1$.}
\end{figure}

\section{\label{sec:sao} Modeling Small-Amplitude Oscillations}

We now study some of the possible concatenations of small-amplitude oscillation segments as seen in Fig. \ref{fig:sao}.  Tangencies of the vector field with the cross-section are given by curves which partition the section into disconnected subsets. The partition that does not contain the attracting and repelling spirals is mapped with full rank to the remaining partition in one return (Fig. \ref{fig:2dmap}(b)), allowing us to restrict our analysis to an invariant two-dimensional subset where the vector field is transverse everywhere. 

Mixed-rank behavior occurs in this subset, as shown in Figure \ref{fig:2dmap}. Note that this set of figures is plotted at a slightly different parameter set from that in Figure \ref{fig:z0}, the main difference being that the line of vector field tangencies intersects a portion of the attracting spiral. Other features of the dynamics persist, including the fact that trajectories jumping left to $S^{a-}_{\eps}$ reach a greater maximal height ($y$-component) than the trajectories jumping right to $S^{a+}_{\eps}$.  As shown in Figure \ref{fig:2dmap}(a), the region is partitioned according to three criteria: their location with respect to the curve of tangency, and their location with respect to the repelling spiral (corresponding to left or right jumps), and their winding number (defined later in this section). 

We will refine this partition in order to create a symbolic model of the dynamics, but we can already state two significant results: \begin{itemize}
\item {\it Mixed-rank dynamics}. As shown in Fig. \ref{fig:2dmap}(c)-(d), the red and blue regions collapse to $S^{a+}_{\eps}$ within one return. This contraction includes those trajectories that return to the cross-section by first jumping right from $S^r_{\eps}$ to $S^{a+}_{\eps}$ at an intermediate height. Figure \ref{fig:2dmap}(b) shows that the yellow subset returns immediately to this low-rank region. The green subset returns either to the low-rank region or to the yellow region. But note that it does not intersect its image, and furthermore, it intersects the yellow region on a portion of the attracting spiral. Therefore, after at most two returns, the dynamics of the points beginning in $\Sigma_0$ (and which did not map to $B_0$) is characterized by the dynamics on the attracting spiral.
\item {\it Trajectories jumping left or right return differently}. Those trajectories jumping left to $S^{a-}_{\eps}$ return to a tiny segment very close to the center of the $S^{a+}_{\eps}$, as shown in Fig. \ref{fig:2dmap}(d). In contrast, the trajectories jumping right sample the entire spiral of $S^{a+}_{\eps}$, as shown in Fig. \ref{fig:2dmap}(c). Thus, multiple intermediate-height jumps to the right are a necessary ingredient in concatenating small- and medium-amplitude oscillations (arising from right jumps) between large-amplitude excursions (arising from left jumps). 
\end{itemize}   

We now construct a dynamical partition of the cross-section. First we define the winding of a trajectory. Let $s$ and $u$ denote a stable and unstable eigenvector, respectively, of the linearization of the flow at $p_{eq}$. Then consider a cylindrical coordinate system with basis $(u,s,n)$ centered at $p_{eq}$, where $n = u \times s$. The {\it winding} of a given trajectory is the cumulative angular rotation (divided by $2\pi$) of the projection of the trajectory onto the $(u,n)$-plane. The {\it winding number} (or simply {\it number of turns}) of a trajectory is the integer part of the winding.

The cumulative angular rotation depends on both the initial and stopping condition of the trajectory, which in turn depend on the section used. Close to $p_{eq}$, the winding of a trajectory measures winding around $W^s$. This is desirable since most of the rotation occurs as trajectories enter small neighborhoods of $p_{eq}$ along $W^s$. 

\begin{figure}
(a) \includegraphics[width=0.44\textwidth]{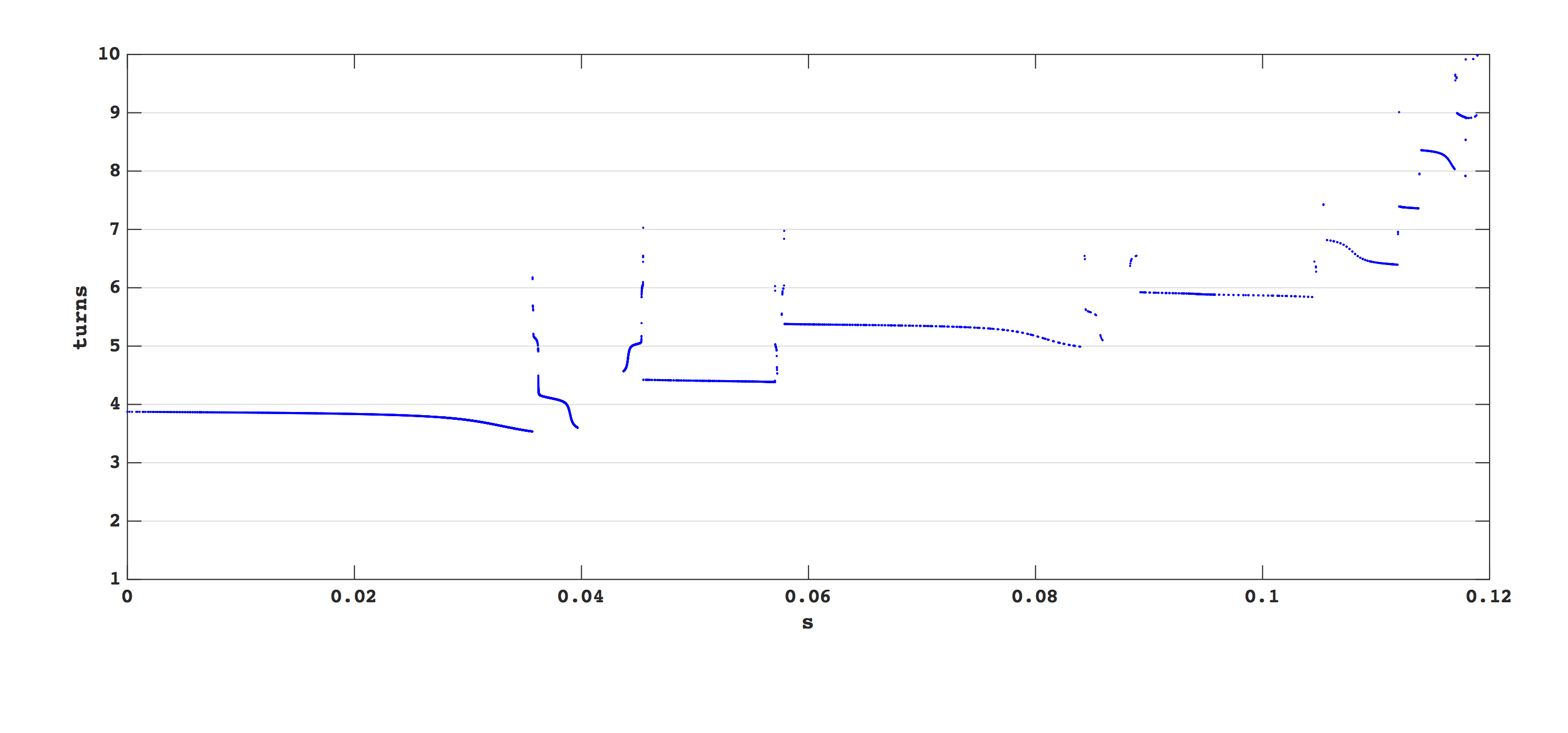}
(b) \includegraphics[width=0.45\textwidth]{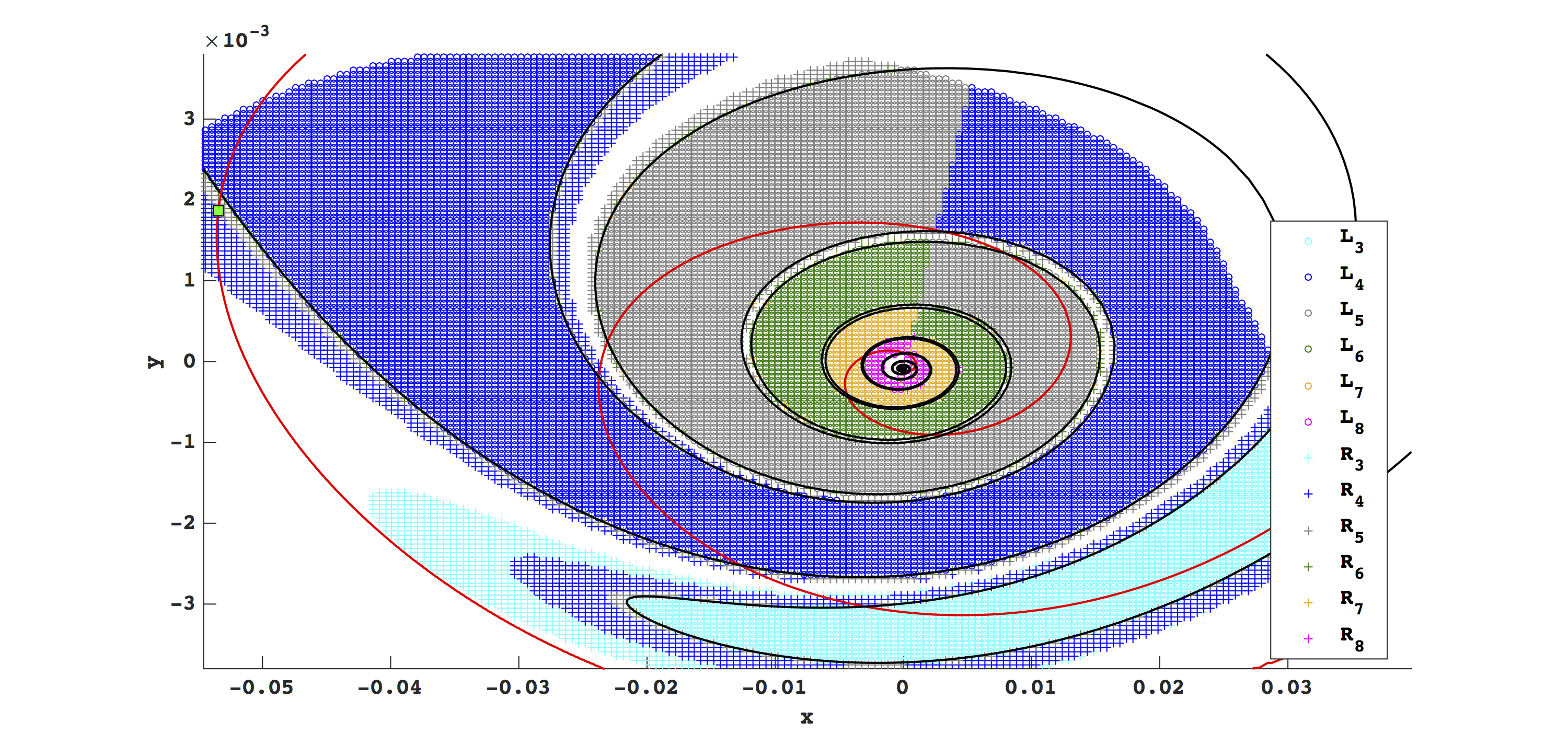}
\caption{\label{fig:turns} (a) Winding of the attracting spiral as a function of its parametrization by arclength. The starting point $s = 0$ is chosen close to the tangency. Positive values of $s$ track the spiral as it turns inward. (b) Partition of the section $\Sigma_0 = \{z = 0\}$ according to number of turns made by corresponding trajectories as well as whether the trajectories turn left or right from $S^r_{\eps}$. Left-turning trajectories are plotted with circles and right-turning trajectories are plotted with crosses. Color definitions: teal, 3 turns; blue, 4 turns; gray, 5 turns; green, 6 turns; gold, 7 turns; magenta, 8 turns. Also plotted are the slow manifolds $S^{a+}_{\eps}$ (red curve) and $S^r_{\eps}$ (black curve) as well as a saddle-point defined in Figure \ref{fig:fporbit} (green square). Parameter set: $\nu \approx 0.00870134$, $a = -0.3$, $b = -1$, $c = 1$.}
\end{figure}

If Figure \ref{fig:turns}(a), we study the winding on a connected subset of the attracting spiral. On this connected subset we may parametrize the spiral by its arclength. The number of turns increases by approximately one whenever $S^{a+}_{\eps}$ intersects $S^{r}_{\eps}$ twice (these intersections occur in pairs since they correspond to bands of trajectories on $S^{a+}_{\eps}$ which leave the region by jumping left to $S^{a-}_{\eps}$). In between these intersections, there are gaps corresponding to regions where $S^{a+}_{\eps}$ intersects $B_0$. 

Sets in the partition are defined according to each trajectory's winding number and jump direction. This partition uses the attracting and repelling spirals as a guide; small rectangles straddling the attracting spiral are contracted strongly transverse to the spiral and stretched along the attracting spiral, giving the dynamics a hyperbolic structure. In the next section we will compute a transverse homoclinic orbit, where this extreme contraction and expansion is shown explicitly.

We restrict ourselves to a subset $S \subset \Sigma_0$ where returns are sufficiently low-rank (i.e. the union of red and blue regions in Fig. \ref{fig:2dmap}(a)). Let $L_{n}\subset S$ (resp. $R_n \subset S$) denote those points whose forward trajectories make $n$ turns before jumping left to $S^{a-}_{\eps}$ (resp. right to $S^{a+}_{\eps}$). Then define $L_{tot} =  \cup_{n=0}^{\infty} L_n$ and $R_{tot} = \cup_{n=0}^{\infty} R_n$. The collection $\mathcal{P} = \{L_i,R_j\}_{i,j=1}^{\infty}$ partitions $S$. 

For a collection of sets $\mathcal{A}$, let $\sigma(\mathcal{A})$ denote the set of all one-sided symbolic sequences $x = x_0 x_1 x_2 \cdots$ with $x_i \in \mathcal{A}$. We can assign to each $x\in S$ a symbolic sequence in $\sigma(\mathcal{P}\cup \{S^c, B_0\})$, also labeled $x$. This sequence is constructed using the return map: $x = \{x_i\}$ is defined by $x_i = \iota(R^i(x))$, where $\iota: \Sigma_0 \to \mathcal{P}\cup\{B_0,S^c\}$ is the natural inclusion map. Note that some symbolic sequences have finite length, as $R$ is undefined over $B_0$. A portion of the partition is depicted in Fig. \ref{fig:turns}.

The results in figures \ref{fig:2dmap} and \ref{fig:turns} and the definition of $B_0$ constrain the allowed symbolic sequences. In particular:

\begin{itemize}
\item there exists a sufficiently large integer $N$ with $R(L_{tot}) \subset S^{a+}_{\eps} \cap (\cup_{n\geq N} L_n \cup R_n \cup B_0)$ (Figs. \ref{fig:2dmap}(d) and \ref{fig:turns}),
\item $R(R_{tot}) \subset S^{a+}_{\eps}\cap \Sigma_0$ (Fig. \ref{fig:2dmap}(c)),
\item $R(S^c) \subset S$ (Fig. \ref{fig:2dmap}(b)), and
\item the set of finite sequences are precisely those containing and ending in $B_0$.
\end{itemize}

The first result implies that for any integer $n \geq 1$, the block $L_n \alpha_m$ (where $\alpha \in \{L,R\}$) is impossible when $m < N$, since $R(L_n)$ is either $B_0$ or $\alpha_{m\geq N}$. For the parameter set used in Fig. \ref{fig:2dmap}, our numerics suggest a lower bound of $N = 13$. On the other hand, the second result reminds us that only right-jumping trajectories are able to sample the entire attracting spiral. The first two results then imply that blocks of type $R_i L_j$ or $R_i R_j$ are necessarily present in the symbolic sequences of orbits which concatenate small-amplitude oscillations with medium-amplitude oscillations as shown in Fig. \ref{fig:sao}, since medium-amplitude oscillations arise precisely from those points on $\Sigma_0$ whose forward trajectories remain bounded away from the saddle-focus (i.e. those points in $\Sigma_0$ sufficiently far from the intersection of $W^s$ with $\Sigma_0$) and jump right. 

The first two results also imply that forward-invariant subsets lie inside the intersection of $S^{a+}_{\eps}$ with $\Sigma_0$. In terms of the full system, it follows that the trajectories corresponding to these points each contain segments which lie  within a sheet of $S^{a+}_{\eps}$. 
	
The attracting spiral and the nonsingular region $S^c$ have a nontrivial intersection. In view of the second result, it is possible for trajectories that we track to sometimes be mapped outside of the subset $S$ where our partition is defined. The third result implies the symbolic sequence of a point $x\in S$ whose forward returns leave the subset $S$ must contain the block
\begin{eqnarray*}
 x_{n_j-1} S^c x_{n_j},
\end{eqnarray*}

where the index $n_j$ is defined by the $j$-th instance when the orbit leaves $S$ and $x_{n_j-1} \in \mathcal{P}$. We can also constrain the possible symbols of $x_{n_j}$ as follows. Trajectories which have intersected $\Sigma_0$ in the singular region $S$ can only return to $S^c$ along the curve $S^{a+}_{\eps}$. Furthermore, our numerical results demonstrate that $R(S^c \cap S^{a+}_{\eps})$ nontrivially intersects subsets of $\mathcal{P}\cup\{B_0\}$ only in the subcollection $\mathcal{P}_c = \{L_3,L_4,L_5,R_3,R_4,R_5,B_0\}$. Therefore $x_{n_j} \in \mathcal{P}_c$ whenever $n_j$ is defined. 

In view of the last result, for each $i \geq 1$ we define the $i$-th {\it escape subset} $E_i$ to be the set of length-$i$ sequences ending in $B_0$. Note that $E_i$ contains the symbol sequences of the points in $R^{-(i-1)}(B_0)$. Section \ref{sec:ret} provides a concrete numerical example of a point in $E_{n}$, where $n$ is at least 1284. 

Let us summarize the main results of the symbolic dynamics. Points which do not have repeating symbolic sequences (i.e. points which are not equilibria or do not belong to cycles) either terminate in $B_0$, indicating that trajectory tends asymptotically to the small-amplitude stable periodic orbit $\Gamma$, or the sequence is infinitely long. In either case, due to the apparent hyperbolicity of the map we can track neighborhoods of points until they leave the set $S$ where the partition is defined. However, those `lost' points on the orbit return to $S$ in one interate and we can resume tracking them.

\section{\label{sec:homorbit2d} Invariant Sets of the Two-Dimensional Return Map}

\begin{figure}
(a)\includegraphics[width=0.46\textwidth]{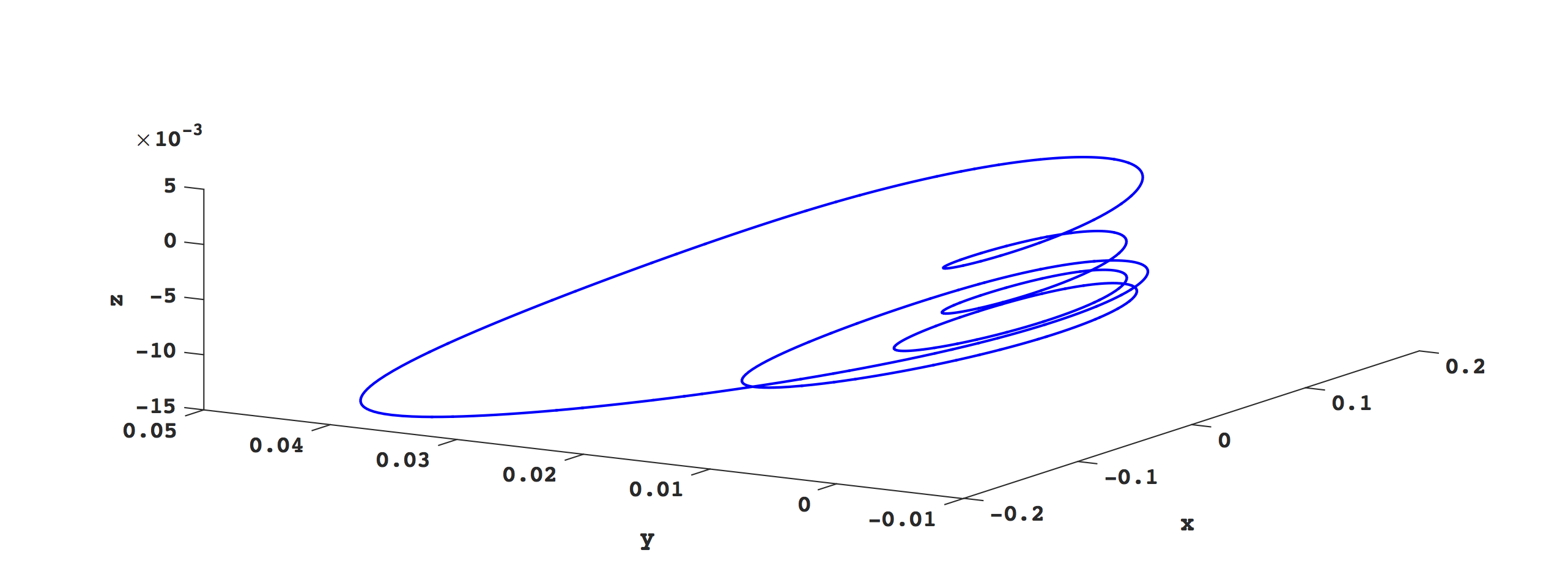}
(b)\includegraphics[width=0.46\textwidth]{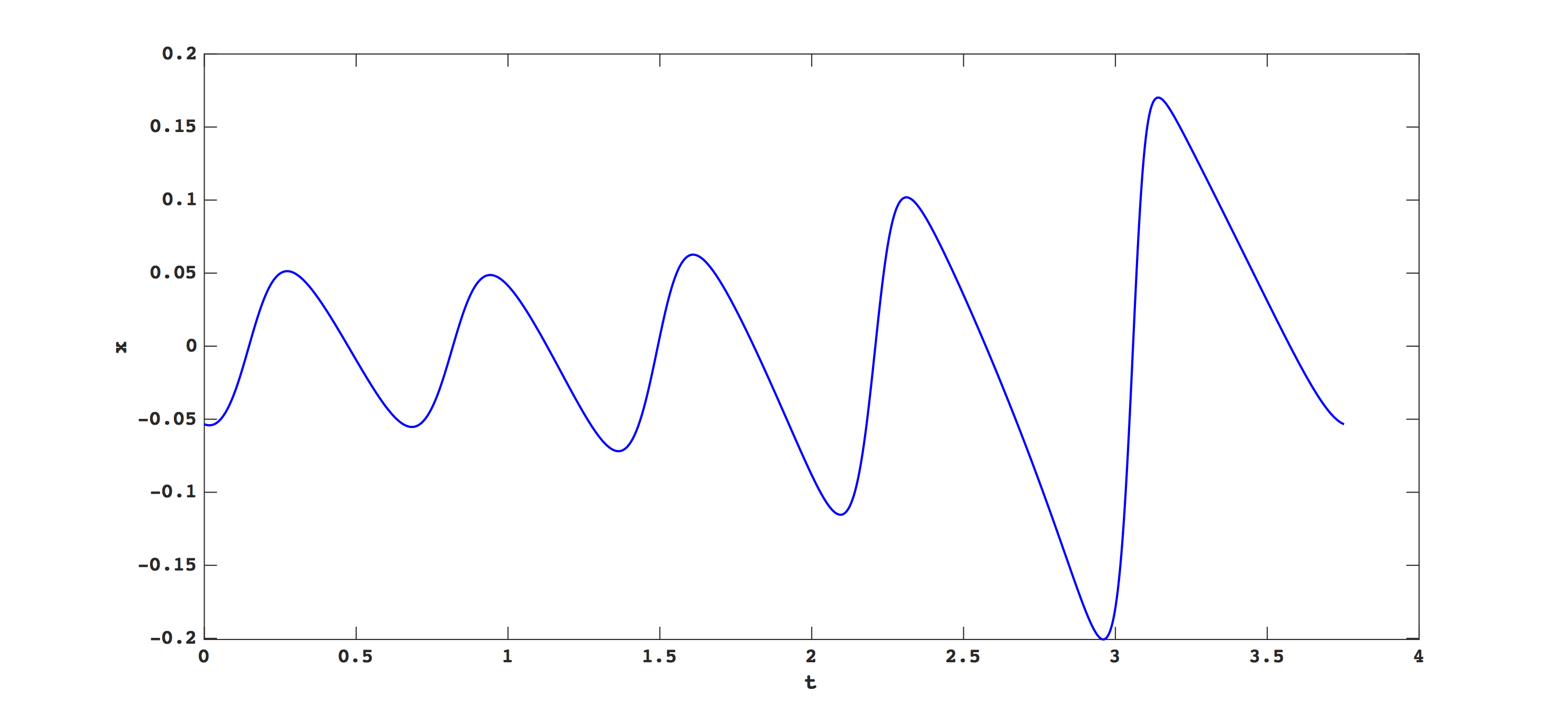}
\caption{\label{fig:fporbit} (a) Mixed-mode oscillation in phase-space corresponding to the saddle point $p \approx (-0.053438, 0.001873)$ of the return map defined on $\Sigma_0 = \{z = 0\}$ and (b) the time series of its $x$-component. Parameter set: $\nu \approx 0.00870134$, $a = -0.3$, $b = -1$, $c = 1$.}
\end{figure}

The structure of the invariant sets of the two-dimensional return map is a key dynamical question, and is related to the intersection of the basin of attraction of the small-amplitude stable periodic orbit with $\Sigma_0$. In this section we focus on two types of invariant sets : fixed points and transverse homoclinic orbits.

\begin{figure}
(a)\includegraphics[width=0.44\textwidth]{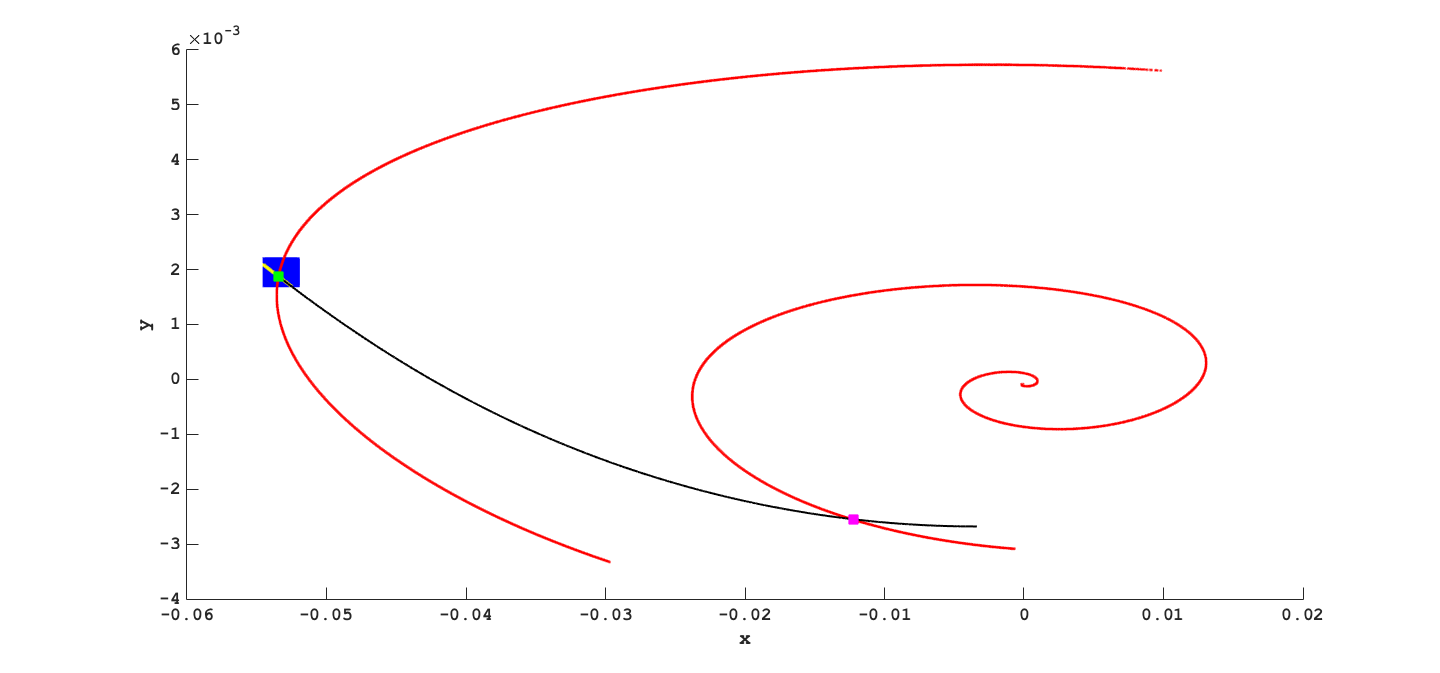}
(b)\includegraphics[width=0.44\textwidth]{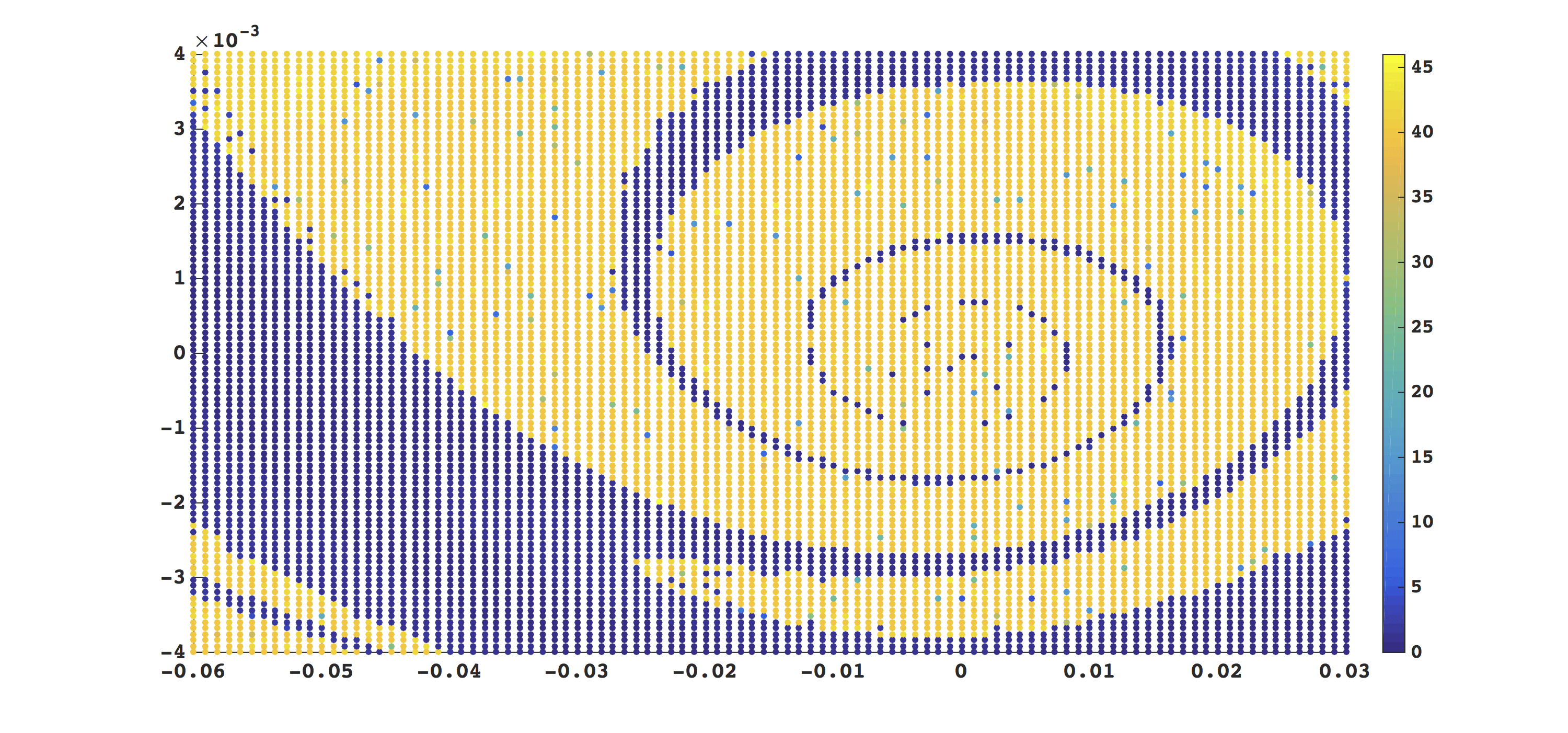}\\
(c)\includegraphics[width=0.44\textwidth]{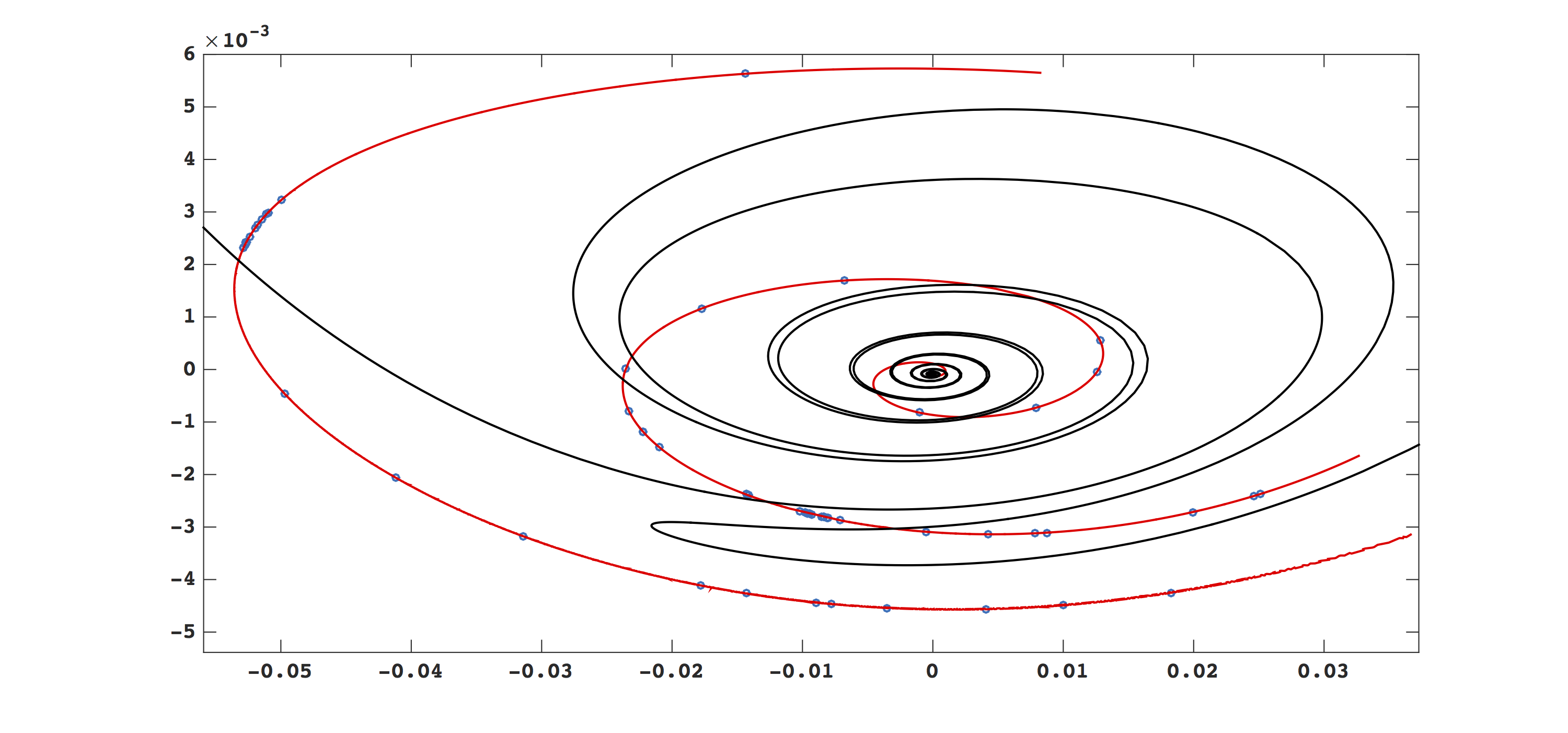}
\caption{\label{fig:homorbit} (a) A saddle equilibrium (green point) of the return map defined on $\Sigma_0 = \{z = 0\}$, together with a neighborhood $U$ (blue grid), image $R(U)$ (red), subset of preimage $U \cap R^{-1}(U)$ (yellow), and a branch of its stable manifold $W^s(p)$(black). The intersection of $R(U)$ with $W^s(p)$ is also shown (magenta point). (b) Color plot of $10^4$ initial conditions beginning in $\Sigma_0$ on a $100\times 100$ grid, whose forward trajectories are integrated for the time interval $t \in [0,600]$. Color denotes number of intersections with $\Sigma_0$ with orientation $\dot{z} < 0$. (c) Last recorded intersection (blue circles) of each trajectory defined in (b) with $\Sigma_0$. The attracting and repelling spirals (red and black curves, respectively) are overlaid. Parameter set: $\nu \approx 0.00870134$, $a = -0.3$, $b = -1$, $c = 1$.}
\end{figure}

Certain invariant sets of the map may be used to construct open sets of points all sharing the same initial block in their symbolic sequence. We briefly describe how the simplest kind of invariant set-- a fixed point-- implies that neighborhoods of points must have identical initial sequences of oscillations. In Figure \ref{fig:fporbit} we plot the saddle-type MMO corresponding to a saddle equilibrium point $p$, whose location in the section $\{z = 0\}$ is plotted in Figures \ref{fig:turns} and \ref{fig:homorbit}. According to Figure \ref{fig:turns}, $p$ has symbolic sequence $R_5 R_5 R_5 \cdots$, in agreement with the time-series shown in Figure \ref{fig:fporbit}(b). We also observe that the fixed-point is sufficiently far away from $W^s$ (the stable manifold of the saddle-focus) that the oscillations remain bounded away from $p_{eq}$. Furthermore, the dynamics in small neighborhoods of $p$ are described by the linearization of the map $R$ near $p$. This implies that small neighborhoods of $p$ consist of points with initial symbolic blocks of $R_5$, where the length of this initial block can be as large as desired. We can relax the condition that this be the initial block by instead considering preimages of these neighborhoods. 

From this case study we observe that arbitrarily long chains of small-amplitude oscillations can be constructed using immediate neighborhoods of fixed points, periodic points, and other invariant sets lying in $S^{a+}_{\eps} \cap \{z = 0\}$. These in turn correspond to complicated invariant sets in the full three-dimensional system. Consequently, the maximum number of oscillations produced by a periodic orbit having one large-amplitude return can be very large at a given parameter value, depending on the number of maximum possible returns to sections in the region containing these local mechanisms. This situation should be compared to earlier studies of folded-nodes, in which trajectories with a given number of small-amplitude oscillations can be classified\cite{wechselberger2005}; and the Shilnikov bifurcation in slow-fast systems, in which trajectories have unbounded numbers of small-amplitude oscillations as they approach the homoclinic orbit.\cite{guckenheimer2015}

In the present system, we observe numerically that the return map $R$ contracts two-dimensional subsets of the cross-section to virtually one-dimensional subsets of $S^{a+}_{\eps}$. Subsequent returns act on $S^{a+}_{\eps}$ by stretching and folding multiple times before further extreme contraction of points transverse to $S^{a+}_{\eps}$. 

If horseshoes exist, they will also appear to break the diffeomorphic structure of the return map due to the strong contraction (although topologically the horseshoe is still given by two-dimensional intersections of forward images of sets with their inverse images under the return map). We can draw a useful analogy to the H\'enon family $H_{a,b}(x,y) = (1-ax^2+y,bx)$ when a strange attractor exists at parameter values $0 < b \ll 1$. The local structure of the strange attractor, outside of neighborhoods of its folds, is given by $C \times \mathbb{R}$, where $C$ is a Cantor set.

Let us now provide a visual demonstration of these issues. Let $U$ be a small neighborhood of the saddle fixed point $p$ that we located in the previous section. In Fig. \ref{fig:homorbit} we plot $U$, $R(U)$, $U\cap R^{-1}(U)$, and $W^s(p)$ on the section $\Sigma_0$. The image is a nearly one-dimensional subset of $S^{a+}_{\eps}$ and the preimage is a thin strip which appears to be foliated by curves tangent to $S^r_{\eps}$. The subsets $R(U)$ and $R^{-1}(U)$ contain portions of $W^u(p)$ and $W^s(p)$, respectively. The transversal intersection of $R(U)$ with $W^s(p)$ is also indicated in this figure.

Numerically approximating the diffeomorphism $R^{-1}$ is a challenging problem. Trajectories which begin on the section and approach the attracting slow manifolds $S^{a\pm}_{\eps}$ in reverse time are strongly separated, analogous to the scenario where pairs of trajectories in forward time are strongly separated by $S^r_{\eps}$. This extreme numerical instability means that trajectories starting on the section and integrated backward in time often become unbounded. In order to compute $W^s(p)$, we therefore resort to a continuation algorithm which instead computes orbits in forward time. We take advantage of the singular behavior of the map to reframe this problem as a boundary value problem, with initial conditions beginning in a line on the section and ending `at' $p$. Beginning with a point $y_0$ along $W^s(p)$, we construct a sequence $\{y_0 , y_1, \cdots\}$ along $W^s(p)$ as follows.

(C1) {\it Prediction step}. Let $w_i = y_{i-1} + h v_i$, where $h$ is a fixed step-size and $v_i$ is a numerically approximated tangent vector to $W^s(p)$ at $y_{i-1}$.

(C2) {\it Correction step}. Construct a line segment $L_i$ of initial conditions perpendicular to $v_i$. Use a bisection method to locate a point $y_i \in L_i$ such that $|R(y_i) - p| < \eps$, where $\eps$ is a prespecified tolerance.

The relevant branch of $W^s(p)$ which intersects $R(U)$ lies inside the nearly singular region of the return map, so the segment $L_i$ can be chosen small enough that $R(L_i)$ is, to double-precision accuracy, a segment of $S^{a+}_{\eps}$ which straddles $p$. This justifies our correction step above.

It is usually not sufficient to assert the existence of a transverse homoclinic orbit from the intersection of the image sets. But in the present case, these structures are organized by the slow manifolds $S^{a+}_{\eps}$ and $S^r_{\eps}$. The strong contraction onto $S^{a+}_{\eps}$ in forward time implies that the discrete orbits comprising $W^u(p)$ must also lie along this slow manifold. The unstable manifold $W^u(p)$ lies inside a member of the $O(\exp(-c/\eps))$-close family which comprises $S^{a+}_{\eps}$, so the forward images serve as good proxies for subsets of $W^u(p)$ itself. On the other hand, when $U$ is sufficiently small, its preimage $R^{-1}(U)$ appears to be foliated by a family of curves tangent to $S^r_{\eps}$, such that one of the curves contains $W^s(p)$ itself.

The Smale-Birkhoff homoclinic theorem\cite{birkhoff1950,smale1965} then implies that there exists a hyperbolic invariant subset on which the dynamics is conjugate to a subshift of finite type. Note that while we expect fixed points to lie in $S^{a+}_{\eps}$ due to strong contraction, we do not expect fixed points to also lie in $S^r_{\eps}$ in the generic case. We end this result by commenting on its apparent degeneracy of the two-dimensional sets $U,R(U),$ and $R^{-1}(U)$. A classical proof of the Smale-Birkhoff theorem uses the set $V = R^k(U) \cap R^{-m}(U)$ (where $k,m\geq 0$ are chosen such that $V$ is nonempty) as the basis for constructing the Markov partition on which the shift is defined.\cite{guckenheimer1983} Here, $V$ is well-approximated by a curve segment. 

Does a positive-measure set of initial conditions approach this hyperbolic invariant set? In other words, is the set an attractor? While it is difficult to assess the invariance of open sets with finite-time computations, our numerics support the conjecture that most initial conditions lying outside $B_0$ tend to the chaotic invariant set without tending asymptotically to $\Gamma$. In Figs. \ref{fig:homorbit}(b)-(c), we study the eventual fates of a grid of initial conditions beginning on $\Sigma_0$. Fig. \ref{fig:homorbit}(b) shows that even after a relatively long integration time of $t = 600$, most initial conditions in $B_0^c$ are able to return repeatedly to $\Sigma_0$. However, it may simply be that the measure of $(R^{-i} (B_0))^c$ decays extremely slowly to $0$ as $i$ tends to infinity. In Fig. \ref{fig:homorbit}(c), we plot the last recorded intersection with $\Sigma_0$ of those trajectories that do not tend to $\Gamma$ within $t = 600$. Even with a relatively sparse set of $10^4$ points, we observe that these intersections sample much of the attracting spiral. Many of the points are not visible at the scale of the figure because they sample the segment shown in Fig. \ref{fig:2dmap}(d) (i.e. the penultimate intersections resulted in the trajectory jumping left to $S^{a-}_{\eps}$).

We now turn to a section transverse to $S^{a+}_{\eps}$ and approximate returns to this section by a one-dimensional map. The advantage of this low-dimensional approximation is that we can readily identify classical bifurcations and routes to complex behavior. We may interpret invariant sets of these one-dimensional returns (such as fixed points, periodic orbits, and more complicated sets) as large-amplitude portions of the mixed-mode oscillations in the full system.

\begin{figure}
(a) \includegraphics[width=0.44\textwidth]{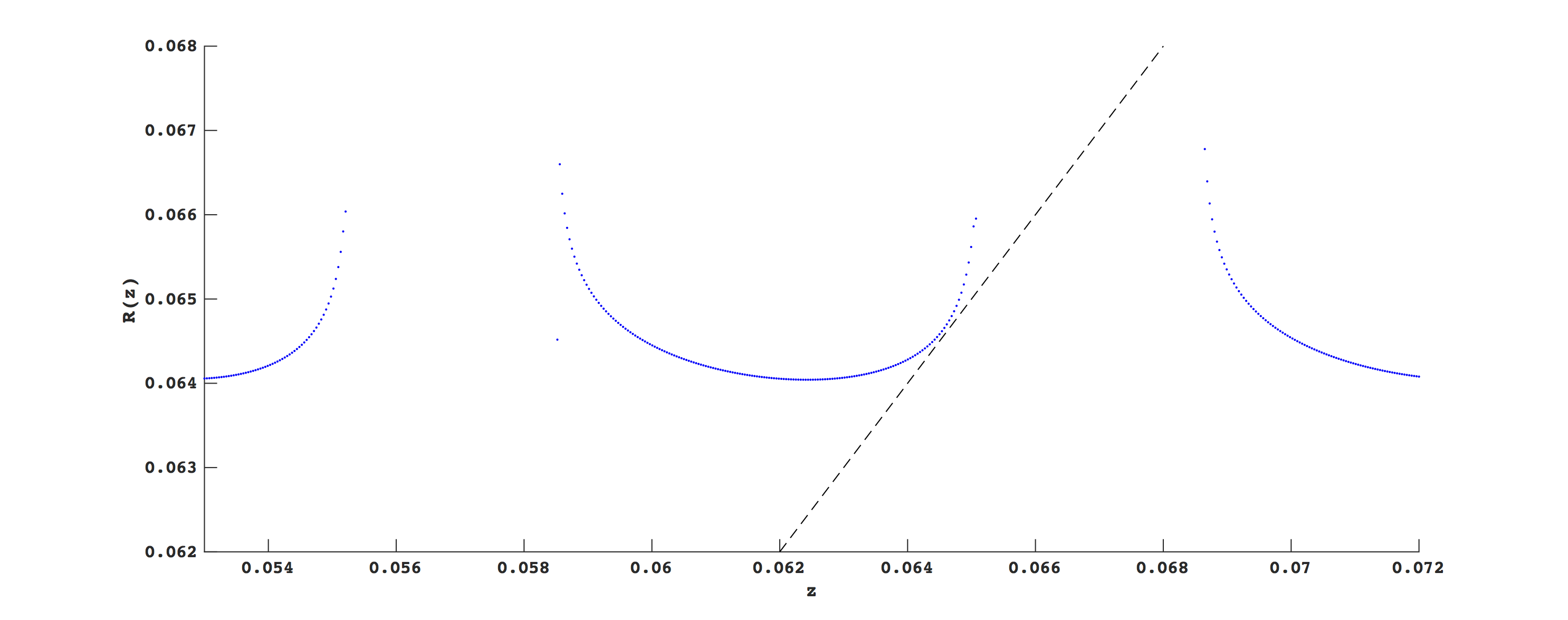}
(b) \includegraphics[width=0.44\textwidth]{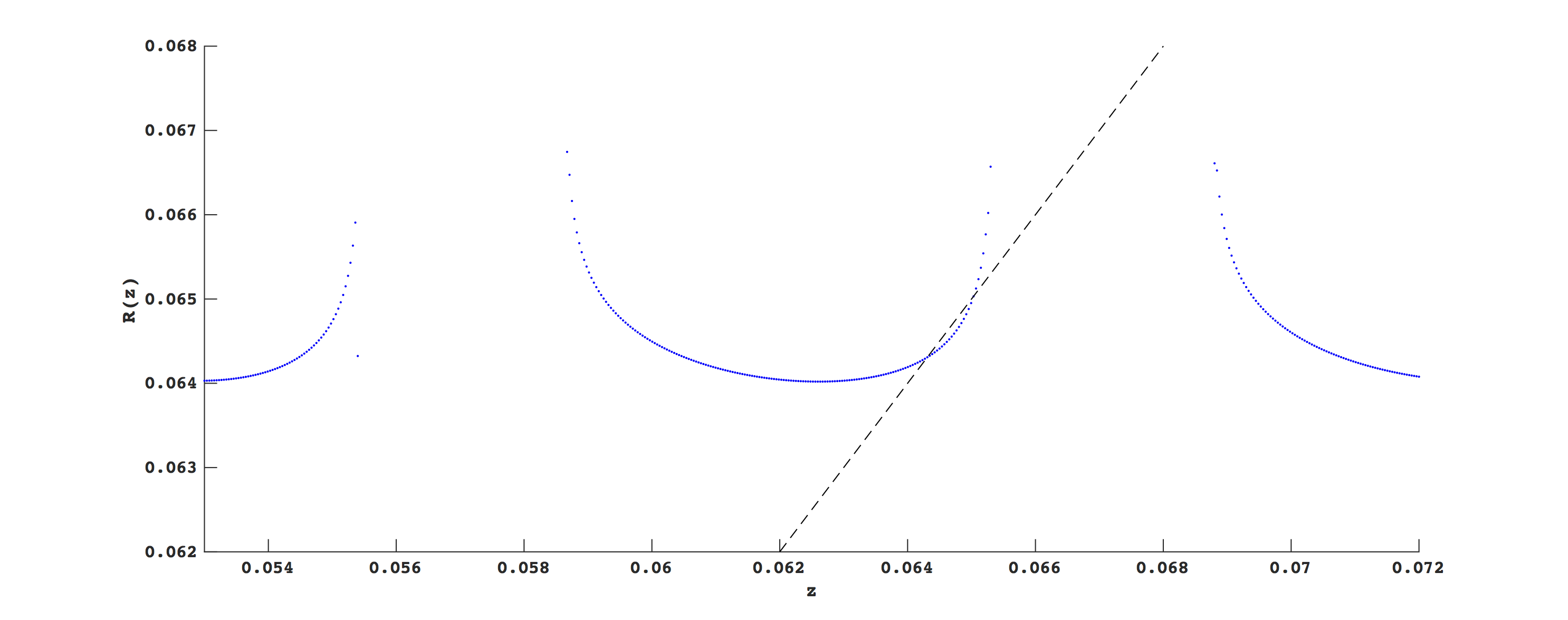}
\caption{\label{fig:sn} Saddle-node bifurcation of periodic orbits in system \eqref{eq:shnf}. (a) $\nu = 0.00801$, (b) $\nu = 0.00802$.  Dashed black line denotes line of fixed points $R(z) = z$. Remaining parameters: $a = -0.3, b = -1, c = 1$.}
\end{figure}

\section{\label{sec:bif} Bifurcations of the One-Dimensional Return Map}

Fixed points of a return map defined on the section $\Sigma_+ = \{x = 0.3\}$ are interpreted in the full system as the locations of mixed-mode oscillations, formed from trajectories making one large-amplitude passage after interacting with the local mechanisms near $L_0$. Similarly, periodic orbits of the (discrete) return map can be used to identify mixed-mode oscillations having more than one large-amplitude passage. We demonstrate common bifurcations associated with these invariant objects. First we locate a saddle-node bifurcation of periodic orbits, in which a pair of orbits coalesce and annihilate each other at a parameter value. 

Figure \ref{fig:sn} demonstrates the existence of a fixed point $z = R(z)$ with unit derivative as $\nu$ is varied within the interval $\left[0.00801, 0.00802\right]$ (remaining parameters are as in Figure \ref{fig:retmap}). This parameter set lies on a generically codimension one branch in the parameter space. We also note that saddle-node bifurcations serve as a mechanism to produce stable cycles in the full system, which in turn may undergo torus bifurcations and period-doubling cascades as a parameter is varied.

\begin{figure}
(a) \includegraphics[width=0.43\textwidth]{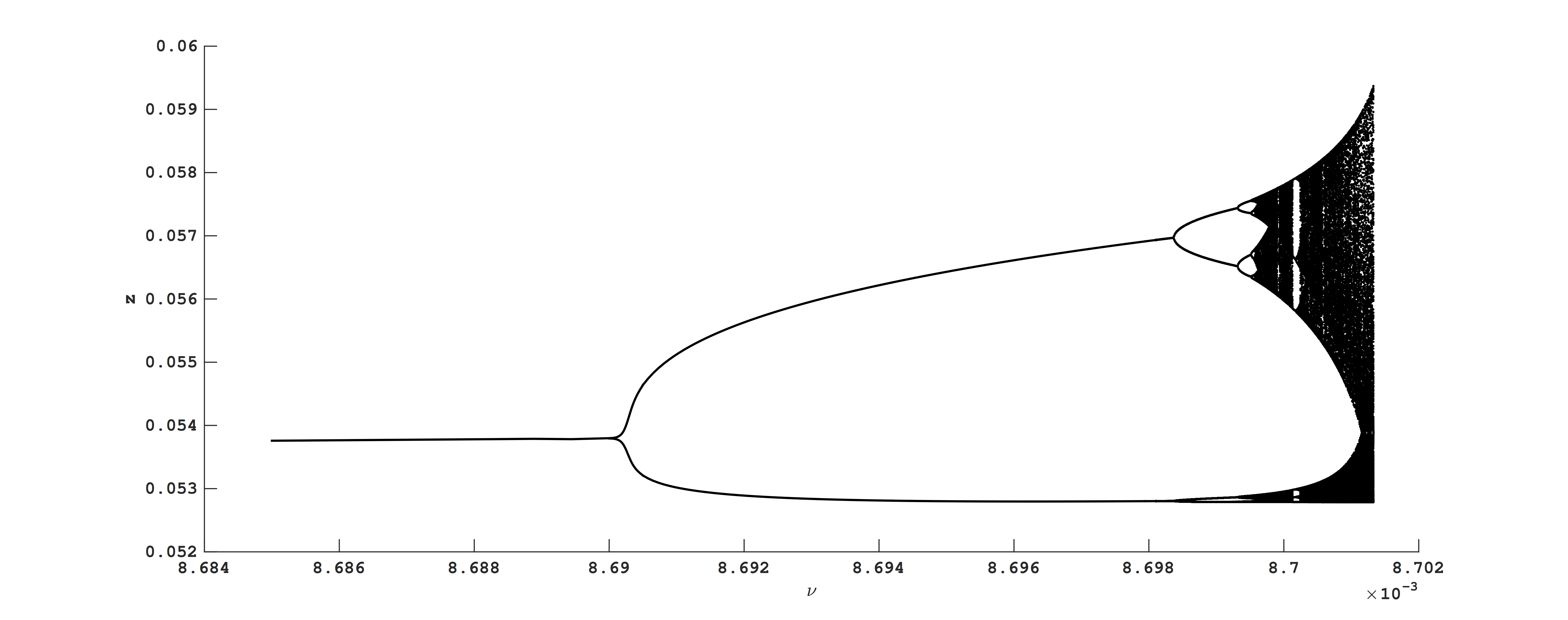}
(b) \includegraphics[width=0.43\textwidth]{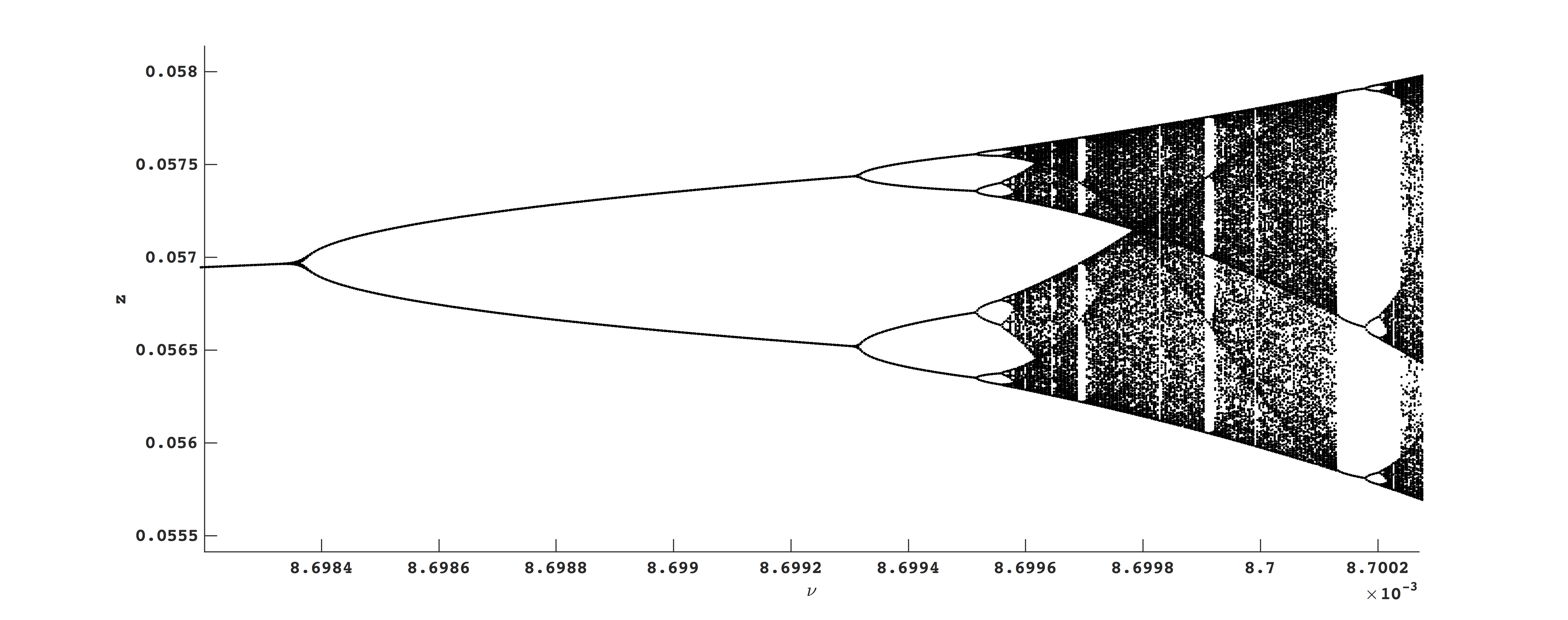}
\caption{\label{fig:bifdiag} (a) Bifurcation sequence of the one-dimensional approximation of the return map $R: \Sigma_+ \to \Sigma_+$ as the parameter $\nu$ is varied from $0.008685$ to $0.0087013$. Remaining parameters are as in Figure \ref{fig:tangbif}. (b) Magnification of upper branch of first period doubling cascade. }
\end{figure}

\begin{figure}[b]
(a) \includegraphics[width=0.42\textwidth]{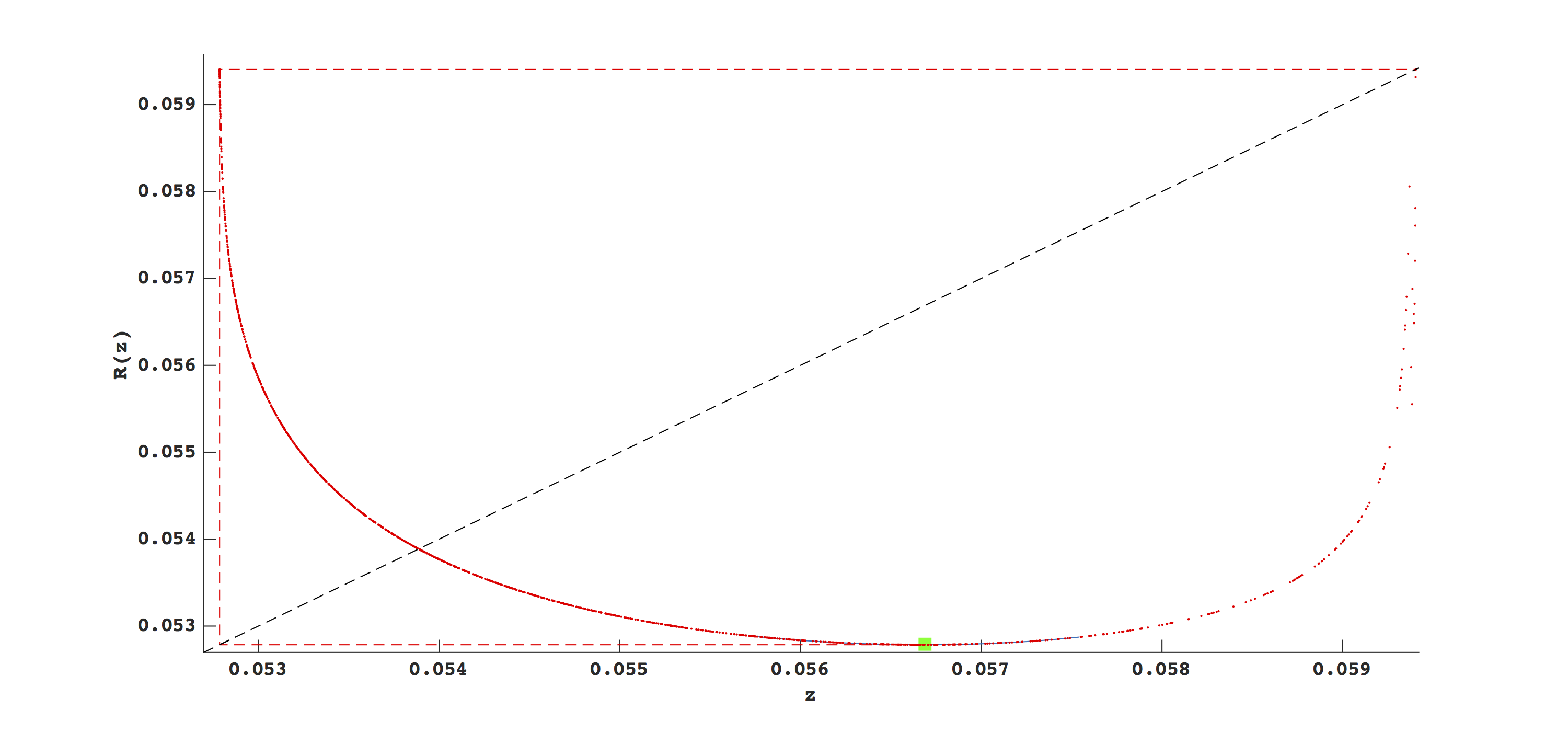}
(b) \includegraphics[width=0.45\textwidth]{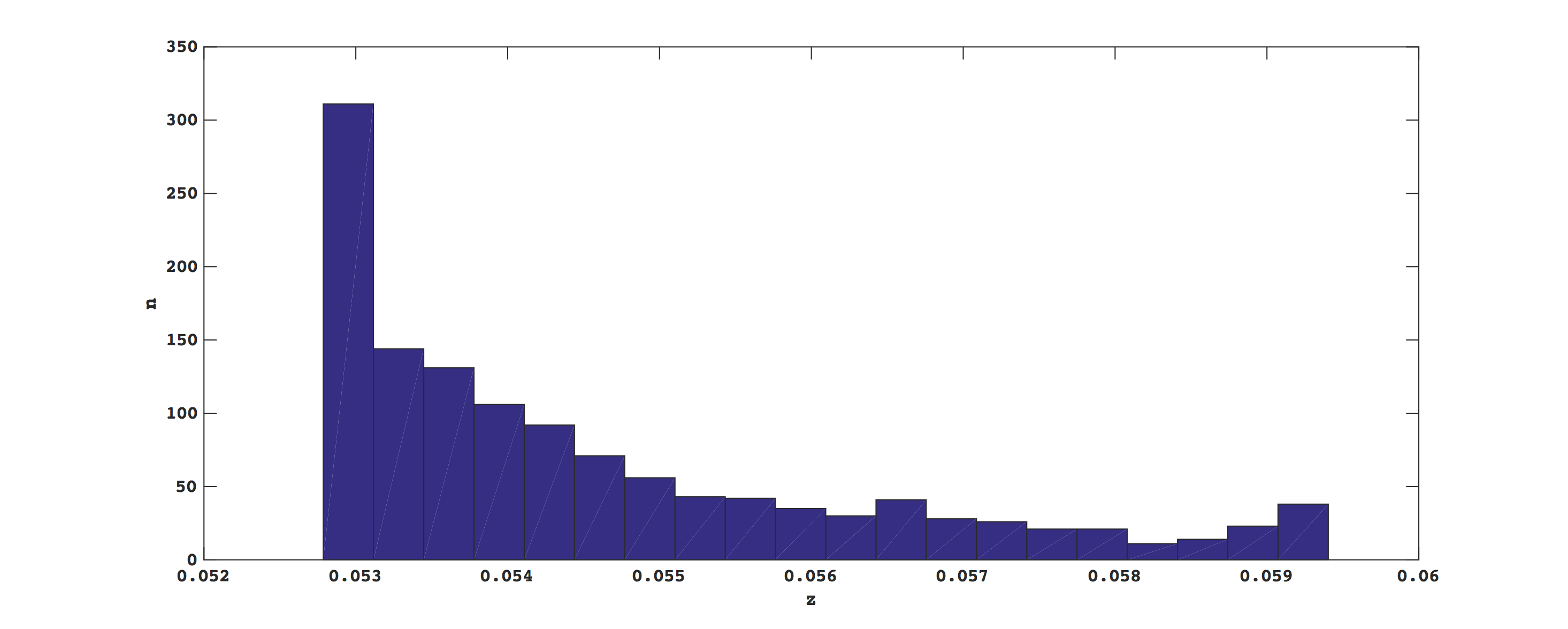}
\caption{\label{fig:crit}  (a) Forward trajectory (red points) of the critical point (green square) under the return map $R: \Sigma_+ \to \Sigma_+$. Red dashed lines indicate the cobweb diagram of the first two iterates of the trajectory to guide the eye. Black dashed line is the line of fixed points. All 1284 forward iterates are plotted. The subsequent iterate lands outside the domain of $R$: the corresponding portion of the full trajectory of \eqref{eq:shnf}  tends asymptotically to $\Gamma$ without returning to $\Sigma_+$. (b) Distribution of points in the forward orbit of the critical point. Parameter set as in Figure \ref{fig:retmap}.}
\end{figure}

The beginning of a period-doubling cascade is identified in the return map $R$ as $\nu$ is varied in the interval $\left[ 0.008685, 0.0087013\right]$ (Figure \ref{fig:bifdiag}a). Within this range, period-3, period-5, and period-6 parameter windows are readily identifiable in Figure \ref{fig:bifdiag}b. The local unimodality of the return map suggests that our ($\nu$-parametrized) family of return maps share some universal properties with maps of the interval that exhibit period-doubling cascades,\cite{feigenbaum1978,coullet1978} despite the nonlinearity at the right boundary of the interval observed in Figure \ref{fig:retmap2}. The cascading structure is clearly robust to small boundary perturbations of the quadratic-like maps we consider. 

We stress that these bifurcations produce additional {\it large-amplitude} oscillations of MMOs. As the parameter is varied between period-doubling events, more small-amplitude twists may be generated. The connection between Figs. \ref{fig:bifdiag} and \ref{fig:turns} is that between each large-amplitude passage, the number and type of small-amplitude twists is determined by the location of periodic points of the return map $R:\Sigma_0\to\Sigma_0$. 

\section{\label{sec:ret} Returns of the critical point}

We recall a classical result of unimodal dynamics for the quadratic family $f_a(x) = 1 - ax^2$ near the critical parameter $a = 2$, where $f_a: I \to I$ is defined on its invariant interval $I$  (when $a = 2$, $I =  \left[ -1,1\right]$). On positive measure sets of parameters near $a = 2$, the map $f_a$ admits absolutely continuous invariant measures with respect to Lebesgue measure.\cite{jakobson1981} These facts depend on the delicate interplay between stretching behavior away from neighborhoods of the critical point, together with recurrence to the arbitrarily small neighborhoods of the critical point as trajectories are `folded back' by the action of $f$. This motivates our current objective: to locate a parameter set for which (i) there exists a forward-invariant subset $\Sigma_u \subset \Sigma_+ $ where $R: \Sigma_u \to \Sigma_u$ has exactly one critical point $c \in \Sigma_u$, and (ii) $R^2(c)$ is a fixed point of $R$.

We couldn't locate a parameter set satisfying both (i) and (ii), but we can obtain a parameter set where $R$ has the topology of Figure \ref{fig:retmap2} (i.e. is unimodal over a sufficiently large interval) and admits a critical point satisfying (ii). This parameter set is numerically approximated using a two-step bisection algorithm. First, a bisection method is used to approximate the critical point $c$ by refining the region where $R'$ first changes sign up to a fixed error term, which we take to be $10^{-15}$. Another bisection method is used to approximate the parameter value at which $|R^2(c) - R^3(c)|$ is minimized. We were able to minimize this distance to $2.5603\times 10^{-8}$ at the parameter value $\nu = 0.0087013381084$, where the remaining parameters are given in Figure \ref{fig:retmap}.

Figure \ref{fig:crit}a depicts the forward trajectory of the critical point near the line of fixed points at this parameter value. The itinerary of $c$ is finite, eventually landing in a subinterval of $\Sigma_u$ where $R$ is undefined. Even so, its forward orbit is unpredictable and samples the interval $\left[ R(c), R^2(c)\right]$ with a nontrivial `transient' density for 1284 iterates (Figure \ref{fig:crit}b). The length of the itinerary is extremely sensitive to tiny ($O(10^{-14})$) perturbations of the parameter $b$, reflecting the sensitive dependence of initial conditions in the selected parameter neighborhood. However, the normalized distributions of the forward iterates behave much more regularly: they are all similar to that shown in Figure \ref{fig:crit}b. We conjecture that there exists a path in parameter space such that these distributions rigorously converge to a smooth distribution on the subset where $R$ is defined.

\section{Concluding remarks}

We have classified much of the complex dynamics arising from a tangency of a slow manifold with an unstable manifold of an equilibrium point. The key to this analysis has been the identification of global bifurcations in carefully-chosen return maps of the system. In particular, transverse homoclinic orbits and period-doubling cascades are identified as mechanisms leading to chaotic behavior in the present system. Our objective has not been to attempt rigorous proofs of the results. Instead, we show that these bifurcations can be identified with fairly standard numerical integration and bisection procedures. The challenge, which is typical in studies of systems with a strong timescale separation, is to use techniques which bypass the numerical instability that occurs in integrating in forward or reverse time.  This is particularly relevant in the continuation procedure that is used to locate transverse homoclinic orbits for the two-dimensional return map, as well as to identify a section high up on $S^{a+}_{\eps}$ which admits an approximation by a one-dimensional map. 

We also motivate the study of maps having the topology shown in Figure \ref{fig:retmap}(a)-(b). These maps are distinguished by two significant features: they admit small disjoint escape subsets, and they are unimodal over most---but not all---of the remainder of the subset over which the map is defined. Sections \ref{sec:bif} and \ref{sec:ret} can then be regarded retrospectively as an introduction to the dynamics of these maps, especially as they compare to the dynamics of unimodal maps. In particular, we observe that such maps undergo period-doubling cascades (Figure \ref{fig:bifdiag}) as a system parameter is varied. The forward trajectory of the critical point is also seen to have a transient density for a range of parameters (for eg., Figure \ref{fig:crit}). These results lead us to conjecture whether absolutely continuous invariant measures and universal bifurcations for unimodal maps persist weakly for the family of maps studied in this paper. The geometric theory of rank-one maps pioneered by Wang and Young\cite{wang2008} is a possible starting point to prove theorems in this direction. This theory has been used successfully to identify chaotic attractors in families of slow-fast vector fields with one fast and two slow variables.\cite{guckenheimer2006} Their technique is based upon approximating returns by one-dimensional maps. 

\begin{acknowledgments}
This work was supported by the National Science Foundation (Grant No. 1006272). The author thanks John Guckenheimer for useful discussions.
\end{acknowledgments}

\bibliography{shnfsaa4-preprint3.bib}

\end{document}